\documentclass[aos,seceqn,dvips]{arximspdf}
\usepackage{graphics}
\doi{10.1214/07-AOS539}
\volume{36}
\issue{5}
\pubyear{2008}
\firstpage{2153}
\lastpage{2182}

\makeatletter
\newproclaim{rmk}{Remark}[section]
\newtheorem{coro}{Corollary}[section]
\newtheorem{theo}{Theorem}[section]
\newtheorem{lemma}{Lemma}

\def\SL{\mathcal{L}}
\def\ra{\rightarrow}
\def\rap{\stackrel{\mathrm{p}}{\ra}}
\def\ms{\sigma}
\def\TT{^\top}
\def\Ex{\mathbb{E}}
\def\SN{\mathcal{N}}
\def\raL2{\stackrel{L_2}{\ra}}
\def\SB{\mathcal{B}}
\def\SC{\mathcal{C}}
\def\SD{\mathcal{D}}
\def\SF{\mathcal{F}}
\def\SO{\mathcal{O}}
\def\SS{\mathcal{S}}
\def\SX{\mathcal{X}}
\def\var{\operatorname{var}}
\def\ma{\alpha}
\def\ml{\lambda}
\def\mg{\gamma}
\def\mp{\partial}

\makeatother

\begin{document}
\begin{frontmatter}

\title{A class of R\'enyi information estimators for multidimensional densities}
\runtitle{A class of R\'enyi information estimators}

\begin{aug}
\author[A]{\fnms{Nikolai} \snm{Leonenko}\ead[label=e1]{leonenkon@cardiff.ac.uk}},\thanksref{aut1}
\thankstext{aut1}{Supported by EPSRC Grant RCMT119.}
\author[B]{\fnms{Luc} \snm{Pronzato}\corref{}
\ead[label=e2]{pronzato@i3s.unice.fr}\thanksref{aut2}}
\thankstext{aut2}{Supported in part by the IST Programme of the European Community under
the PASCAL Network of Excellence, IST-2002-506778.}
and
\author[A]{\fnms{Vippal} \snm{Savani}\ead[label=e3]{savaniv@cardiff.ac.uk}}
\runauthor{N. Leonenko, L. Pronzato and V. Savani}
\affiliation{Cardiff University,
CNRS/Universit\'e de Nice--Sophia Antipolis and~Cardiff~University}
\address[A]{N. Leonenko\\
V. Savani \\
Cardiff School of Mathematics \\
Cardiff University \\
Senghennydd Road \\
Cardiff CF24 4AG \\
United Kingdom\\
\printead{e1}\\
\phantom{E-mail: }\printead*{e3}} 
\address[B]{L. Pronzato\\
Laboratoire I3S\\
Les Algorithmes\\
CNRS/Universit\'e de Nice--Sophia Antipolis\\
2000 route des Lucioles\\
BP 121\\
06903 Sophia-Antipolis Cedex\\
France\\
\printead{e2}}
\end{aug}

\received{\smonth{7} \syear{2007}}
\revised{\smonth{7} \syear{2007}}

%
\begin{abstract}
A class of estimators of the R\'enyi and Tsallis entropies of an
unknown distribution $f$ in $\mathbb{R}^m$ is presented. These
estimators are based on the $k$th nearest-neighbor distances
computed from a sample of $N$ i.i.d. vectors with distribution $f$.
We show that entropies of any order $q$, including Shannon's
entropy, can be estimated consistently with minimal assumptions on
$f$. Moreover, we show that it is straightforward to extend the
nearest-neighbor method to estimate the statistical distance between
two distributions using one i.i.d.\ sample from each.
\end{abstract}

%


\end{frontmatter}

\section{Introduction}

We consider the problem of estimating the R\'enyi \cite{Renyi61} entropy
%
\begin{equation}\label{Renyi}
H^*_q = \frac{1}{1-q} \log\int_{\mathbb{R}^m} f^q(x)
\,dx , \qquad q\neq1 ,
\end{equation}
or the Havrda and Charv\'at \cite{HavrdaC67} entropy (also called Tsallis
\cite{Tsallis88} entropy)
%
\begin{equation}\label{Tsallis}
H_q = \frac{1}{q-1} \biggl( 1 - \int_{\mathbb{R}^m} f^q(x) \,dx
\biggr)
, \qquad q\neq1 ,
\end{equation}
of a random vector $X\in\mathbb{R}^m$ with probability measure $\mu$
which has density $f$ with respect to the Lebesgue measure, from $N$
independent and identically distributed (i.i.d.) samples
$X_1,\ldots,X_N$, $N\geq2$. Note that $H_q^*$ can be expressed as a
function of $H_q$. Indeed, $H_q^* = \log[1-(q-1)H_q]/(1-q)$, and for
any $q$, $d(H_q^*)/d(H_q)>0$ and $[d^2(H_q^*)/d(H_q)^2]/(q-1)>0$.
For $q<1$ and $q>1$, $H_q^*$ is thus a strictly increasing concave
and convex function of $H_q$ respectively and the maximization of
$H_q^*$ and $H_q$ are equivalent. Hence, in what follows we shall
speak indifferently of $q$-entropy maximizing distributions. When
$q$ tends to 1, both $H_q$ and $H^*_q$ tend to the
(Boltzmann--Gibbs-) Shannon entropy
%
\begin{equation}\label{Shannon}
H_1 = - \int_{\mathbb{R}^m} f(x) \log f(x) \, dx .
\end{equation}
We consider a new class of estimators of $H_q$ and $H^*_q$ based on
the approach proposed by
Kozachenko and Leonenko
\cite{KozachenkoL87} who consider the estimation of $H_1$; see
also \cite{GoriaLMNI04}. Within the classification made in
\cite{BeirlantDGvdM97}, which also contains an outstanding
overview of nonparametric Shannon entropy estimation, the method
falls in the category of nearest-neighbor distances. See also
\cite{HallM93}. When $m=1$, the nearest-neighbor method is
related to sample-spacing methods; see, for example, \cite{Vasicek76} for
an early reference concerning Shannon entropy. It also has some
connections with the more recent random-graph approach of
Redmond and Yukich \cite{RedmondY96}, who, on the supposition that the distribution
is supported on $[0,1]^m$ together with some smoothness
assumptions on $f$, construct a strongly consistent estimator of
$H_q^*$ for $0<q<1$ (up to an unknown bias term independent of $f$
and related to the graph properties). For $q\neq1$, our
construction relies on the estimation of the integral
%
\begin{equation}\label{Iq}
I_q=\Ex\{f^{q-1}(X)\} = \int_{\mathbb{R}^m} f^q(x)\, dx
\end{equation}
through the computation of conditional moments of nearest-neighbor
distances. It thus possesses some similarities with that of
Evans, Jones and Schmidt \cite{EvansJS2002}, who establish the weak consistency of an
estimator of $I_q$ for $m \geq2$ and $q<1$ under the conditions
that $f$ is continuous and strictly positive on a compact convex
subset $\SC$ of $\mathbb{R}^m$, with bounded partial derivatives on
$\SC$. In comparison to Redmond and Yukich \cite{RedmondY96} and Evans, Jones and Schmidt \cite{EvansJS2002},
our results cover a larger range of values for $q$ and do not rely
on assumptions of regularity or bounded support for $f$. For the
sake of completeness, we also consider the case $q=1$, that is, the
estimation of Shannon entropy, with results obtained as corollaries
of those for $q\neq1$ (at the expense of requiring slightly
stronger conditions than Kozachenko and Leonenko \cite{KozachenkoL87}).

The entropy (\ref{Tsallis}) is of interest in the study of
nonlinear Fokker--Planck equations, with $q<1$ for the case of
subdiffusion and $q>1$ for superdiffusion; see \cite{TsallisB96}.
Values of $q\in[1,3]$ are used by Alemany and Zanette \cite{AlemanyZ94} to study the
behavior of fractal random walks. Applications for quantizer
design, characterization of time-frequency distributions, image
registration and indexing, texture classification and image matching
etc., are indicated by
Hero et al. \cite{HeroBMG2002}, Hero and Michel~\cite{HeroM99} and Neemuchwala, Hero and Carson \cite{NeemuchwalaHC2005}. Entropy
minimization is used by Pronzato, Thierry and Wolsztynski \cite{PTWmoda72004}, Wolsztynski,  Thierry and Pronzato \cite{WTPa05} for parameter
estimation in semi-parametric models. Entropy estimation is also a
basic tool for independent component analysis in signal
processing; see, for example, \cite{KraskovSG2004,LearnedMillerF2003}.

The entropy $H_q$ is a concave function of the density for $q>0$
(and convex for $q<0$). Hence, $q$-entropy maximizing distributions,
under some specific constraints, are uniquely defined for $q>0$. For
instance, the $q$-entropy maximizing distribution is uniform under
the constraint that the distribution is finitely supported. More
interestingly, for any dimension $m\geq1$, the $q$-entropy
maximizing distribution with a given covariance matrix is of the
multidimensional Student-$t$ type if $m/(m+2)<q<1$; see
\cite{VignatHC2004}. This generalizes the well-known property that
Shannon entropy $H_1$ is maximized for the normal distribution. Such
entropy-maximization properties can be used to derive nonparametric
statistical tests by following the same approach as
Vasicek \cite{Vasicek76} who tests normality with $H_1$; see also
\cite{GoriaLMNI04}.

The layout of the paper is as follows. Section~\ref{S:motivation}
develops some of the motivations and applications just mentioned
(see also Section~\ref{S:divergences} for signal and image
processing applications). The main results of the paper are
presented in Section~\ref{S:main-results}. The paper is focused on
entropy estimation, but in Section~\ref{S:divergences} we show how a
slight modification of the method also allows us to estimate
statistical distances and divergences between two distributions.
Section~\ref{S:examples} gives some examples and
Section~\ref{S:further} indicates some related results and possible
developments. The proofs of the results of
Section~\ref{S:main-results} are collected in Section~\ref{S:Proofs}.

\section{Properties, motivation and applications} \label{S:motivation}
\subsection{Nonlinear Fokker--Planck equation and entropy}
Consider a family of time-dependent p.d.f.'s $f_t$. The p.d.f.
that maximizes R\'enyi entropy (\ref{Renyi}) [and Tsallis entropy
(\ref{Tsallis})] subject to the constraints $\int_\mathbb{R}
f_t(x) \,dx = 1$, $\int_\mathbb{R} [x-\bar x(t)]\times f_t^q(x) \,dx =
0$, $\int_\mathbb{R} [x-\bar x(t)]^2 f_t^q(x) \,dx =
\ms_q^2(t)$, for fixed $q>1$, is the solution of a nonlinear
Fokker-Planck (or Kolmogorov) equation; see \cite{TsallisB96}.

Let $X$ and $Y$ be two independent random vectors respectively in
$\mathbb{R}^{m_1}$ and $\mathbb{R}^{m_2}$. Define $Z=(X,Y)$ and let
$f(x,y)$ denote the joint density for $Z$. Let $f_1(x)$ and $f_2(y)$
be the marginal densities for $X$ and $Y$ respectively, so that
$f(x,y)=f_1(x)f_2(y)$. It is well known that the Shannon and R\'enyi
entropies (\ref{Shannon}) and (\ref{Renyi}) satisfy the additive
property $H_q^*(f)=H_q^*(f_1)+H_q^*(f_2)$, $q\in\mathbb{R}$, while
for the Tsallis entropy (\ref{Tsallis}), one has
$H_q(f)=H_q(f_1)+H_q(f_2)+(1-q)H_q(f_1)H_q(f_2)$. The first property
is known in physical literature as the extensivity property of
Shannon and R\'enyi entropies, while the second is known as
nonextensivity (with $q$ the parameter of nonextensivity).
The paper by Frank and Daffertshofer \cite{FrankD2000} presents a survey of results
related to entropies in connection with nonlinear Fokker--Planck
equations and normal or anomalous diffusion processes. In
particular, the so-called Sharma and Mittal entropy $H_{q,s} =
[ 1- (I_q)^{(s-1)/(q-1)}]/(s-1)$, with $q,s>0$,
$q,s\neq1$ and $I_q$ given by (\ref{Iq}), represents a possible
unification of the (nonextensive) Tsallis entropy (\ref{Tsallis})
and (extensive) R\'enyi entropy (\ref{Renyi}). It satisfies
$\lim_{s\ra1} H_{q,s} = H_q^*$, $\lim_{s,q\ra1} H_{q,s} = H_1$,
$H_{q,q}=H_q$ and $\lim_{q\ra1} H_{q,s} = \{ 1 -\exp[-(s-1)H_1]
\}/(s-1)= H_s^G$, $s>0$, $s\neq1$,
where $H_s^G$ is known as Gaussian entropy. Notice that a
consistent estimator of $H_{q,s}$ can be obtained from the
estimator of $I_q$ presented in Section~\ref{S:main-results}.

\subsection{Entropy maximizing distributions}
\label{S:entropymaximizing}

The $m$-dimensional random vector $X=([X]_1,\ldots,[X]_m)\TT$ is
said to have a multidimensional Student distribution
$T(\nu,\Sigma,\mu)$ with mean $\mu\in\mathbb{R}^m$, scaling or
correlation matrix $\Sigma$, covariance matrix $C=\nu\Sigma/(\nu-2)$
and $\nu$ degrees of freedom if its p.d.f.\ is
%
\begin{eqnarray}\label{f-Student}
f_\nu(x) &=& \frac{1}{ (\nu\pi)^{m/2}
}\nonumber\\[-8pt]\\[-8pt]
&&{}\times \frac{\Gamma((m+\nu)/2)} { \Gamma(\nu/2) }
\frac{1}{|\Sigma|^{1/2}
[1+(x-\mu)\TT[\nu\Sigma]^{-1}(x-\mu)]^{(m+\nu)/2}} ,\nonumber
\end{eqnarray}
$x\in\mathbb{R}^m$. The characteristic function of the distribution
$T(\nu,\Sigma,\mu)$ is
\[
\phi(\zeta)=\Ex\exp(i\langle\zeta,X\rangle) =
\exp(i\langle\zeta,\mu\rangle)
K_{\nu/2}\bigl(\sqrt{\nu\zeta\TT\Sigma\zeta}\bigr)
\bigl(\sqrt{\nu\zeta\TT\Sigma\zeta}\bigr)^{\nu/2}
\frac{2^{1-\nu/2}}{\Gamma(\nu/2)} ,
\]
$\zeta\in\mathbb{R}^m$, where $K_{\nu/2}$ denotes the modified
Bessel function of the second order. If $\nu=1$, then
(\ref{f-Student}) is the $m$-variate Cauchy distribution. If
$(\nu+m)/2$ is an integer, then (\ref{f-Student}) is the $m$-variate
Pearson type VII distribution. If $Y$ is $\SN(0,\Sigma)$ and if $\nu
S^2$ is independent of $Y$ and $\SX^2$-distributed with $\nu$
degrees of freedom, then $X=Y/S+\mu$ has the p.d.f.\
(\ref{f-Student}). The limiting form of (\ref{f-Student}) as
$\nu\ra\infty$ is the $m$-variate normal distribution
$\SN(\mu,\Sigma)$. The R\'enyi entropy (\ref{Renyi}) of
(\ref{f-Student}) is
\begin{eqnarray*}
H_q^* &=& \frac{1}{1-q} \log
\frac{B(q(m+\nu)/2-m/2,m/2)}{B^q(\nu/2,m/2)}\\
&&{}+ \frac12 \log[(\pi\nu)^m|\Sigma|] - \log
\Gamma\biggl(\frac{m}{2}\biggr) ,\qquad q>\dfrac{m}{m+\nu}.
\end{eqnarray*}
It converges as $\nu\ra\infty$ to the R\'enyi entropy
%
\begin{eqnarray}\label{Hqnormal}
H_q^*(\mu,\Sigma) &=& \log[(2\pi)^{m/2}|\Sigma|^{1/2}] -
\frac{m}{2(1-q)} \log q \nonumber\\[-8pt]\\[-8pt]
&=& H_1(\mu,\Sigma) - \frac{m}{2} \biggl(1+ \frac{\log q}{1-q}
\biggr) \nonumber
\end{eqnarray}
of the multidimensional normal distribution $\SN(\mu,\Sigma)$. When
$q\ra1$, $H_q^*(\mu,\Sigma)$ tends to $H_1(\mu,\Sigma)=\log[(2\pi
e)^{m/2}|\Sigma|^{1/2}]$, the Shannon entropy of $\SN(\mu,\Sigma)$.
For $m/(m+2) < q < 1$, the $q$-entropy maximizing distribution under
the constraint
%
\begin{equation}\label{Cconstraint}
\Ex(X-\mu)(X-\mu)\TT=C
\end{equation}
is the Student distribution $T(\nu,(\nu-2)C/\nu,0)$ with
$\nu=2/(1-q)-m>2$; see \cite{VignatHC2004}. For $q>1$, we define
$p=m+2/(q-1)$ and the $q$-entropy maximizing distribution under
the constraint (\ref{Cconstraint}) has then finite support given
by $\Omega_q = \{x\in\mathbb{R}^m \dvtx (x-\mu)\TT[(p+2)C]^{-1}(x-\mu)
\leq1 \}$. Its p.d.f.\ is
%
\begin{eqnarray}\label{f-finite-support}
&&f_p(x)\nonumber\\[4pt]\\[-20pt]
&&\qquad = \cases{
\dfrac{\Gamma(p/2+1)}{|C|^{1/2}[\pi(p+2)]^{m/2}
\Gamma((p-m)/2+1)}  \cr
\qquad{} \times[1-(x-\mu)\TT[(p+2)C]^{-1}(x-\mu)]^{1/(q-1)}, &\quad if
$x\in\Omega_q$
\cr
0, &\quad otherwise.
}\nonumber
\end{eqnarray}
The characteristic function of the p.d.f.\ (\ref{f-finite-support})
is given by
\[
\phi(\zeta) = \exp(i\langle\zeta,\mu\rangle) 2^{p/2}
\Gamma\biggl(\frac{p}{2}+1\biggr) |\zeta\TT(p+2)C\zeta|^{-p/2}
J_{p/2}\bigl(|\zeta\TT(p+2)C\zeta|\bigr) ,
\]
$\zeta\in\mathbb{R}^m$, where $J_{\nu/2}$ denotes the Bessel
function of the first kind.

Alternatively, $f_\nu$ for $q<1$ or $f_p$ for $q>1$ also maximizes the
Shannon entropy (\ref{Shannon}) under a logarithmic constraint; see
\cite{Kapur89,Zografos99}. Indeed, when $q<1$, $f_\nu(x)$ given by
(\ref{f-Student}) with $\nu=2/(1-q)-m$ and $\Sigma=(\nu-1)C/\nu$
maximizes $H_1$ under the constraint
\[
\int_{\mathbb{R}^m} \log\bigl(1+x\TT[(\nu-2)C]^{-1}x\bigr) f(x) \,dx =
\Psi\biggl(\frac{\nu+m}{2}\biggr) - \Psi\biggl(\frac{\nu}{2}\biggr),
\]
and when $q>1$, $f_p(x)$ given by (\ref{f-finite-support}) with
$p=2/(q-1)+m$ maximizes $H_1$ under
\[
\int_{\mathbb{R}^m} \log\bigl(1-x\TT[(p+2)C]^{-1}x\bigr) f(x) \,dx =
\Psi\biggl(\frac{p}{2}\biggr) - \Psi\biggl(\frac{p+m}{2}\biggr) ,
\]
where $\Psi(z)=\Gamma'(z)/\Gamma(z)$ is the digamma function.

\subsection{Information spectrum}
\label{S:spectrum}

Considered as a function of $q$, $H_q^*$ (\ref{Renyi}) is known as
the spectrum of R\'enyi information; see \cite{Song2001}. The
value of $H_q^*$ for $q=2$ corresponds to the negative logarithm
of the well-known efficacy parameter $\Ex f(X)$ that arises in
asymptotic efficiency considerations. Consider now
%
\begin{equation}\label{dotH1def}
\dot H_1 = \lim_{q\ra1} \frac{dH_q^*}{dq} .
\end{equation}
It satisfies
\begin{eqnarray*}
\dot H_1 &=& \lim_{q\ra1} \frac{\log\int_{\mathbb{R}^m}
f^q(x) \,dx}{(1-q)^2} + \frac{ \int_{\mathbb{R}^m} f^q(x) \log
f(x) \,dx}{(1-q)\int_{\mathbb{R}^m} f^q(x) \,dx}  \\
& = & -\frac12 \biggl\{ \int_{\mathbb{R}^m} f(x) [\log f(x)]^2
\, dx - \biggl[ \int_{\mathbb{R}^m} f(x) \log f(x) \,dx \biggr]^2
\biggr\} \\
&=& -\frac12 \var[\log f(X)] .
\end{eqnarray*}
The quantity $S(f)=-2 \dot H_1=\var[\log f(X)]$ gives a measure of
the intrinsic shape of the density $f$; it is a location and scale
invariant positive functional ($S(f) = S(g)$ when $f(x)=\ms^{-1}
g[(x-\mu)/\ms]$). For the multivariate normal distribution
$\SN(\mu,\Sigma)$, $H_q^*$ is given by (\ref{Hqnormal}) and
$S(f)=m/2$. For the one-dimensional Student distribution with $\nu$
degrees of freedom (for which $\Ex X^{\nu-1}$ exists, but not $\Ex
X^{\nu}$), with density
\[
f_\nu(x)= \frac{1}{(\nu\pi)^{1/2} }
\frac{\Gamma(\nu/2+1/2)}{\Gamma(\nu/2)}
\frac{1}{(1+ x^2/\nu)^{(\nu+1)/2}} ,
\]
we obtain
%
\begin{eqnarray}\label{Hstudentm1}
H_q^* &=& \frac{1}{1-q} \log\frac{
B(q(\nu+1)/2-1/2,1/2)}{B^q(\nu/2,1/2)} +
\frac12 \log\nu, \qquad q> \frac{1}{\nu+1} ,
\nonumber\\[-8pt]\\[-8pt]
S({f_\nu}) &=& \cases{
\frac{\pi^2}{3} \simeq3.2899 ,&\quad for $\nu=1$ (Cauchy
distribution), \cr
9- \frac34 \pi^2 \simeq1.5978, & \quad for $\nu=2$, \cr
\frac43 \pi^2 -12 \simeq1.1595, & \quad for  $\nu=3$, \cr
\frac{775}{36} - \frac{25}{12}\pi^2 \simeq0.9661, & \quad for
$\nu=4$,\cr
3\pi^2 - \frac{115}{4} \simeq0.8588, & \quad for $ \nu=5$,
}\nonumber
\end{eqnarray}
and, more generally, $S({f_\nu}) = (1/4) (\nu+1)^2 \{\dot
\Psi(\nu/2) - \dot\Psi[(\nu+1)/2] \}$, with $\dot\Psi(x)$ the
trigamma function, $\dot\Psi(x)=d^2\log\Gamma(x)/dx^2$. The
information provided by $S(f)$ on the shape of the distribution
complements that given by other more classical characteristics
like kurtosis. [Note that the kurtosis is not defined for $f_\nu$
when $\nu\leq4$; the one-dimensional Student distribution $f_6$
and the bi-exponential Laplace distribution $f_L$ have the same
kurtosis but different values of $S(f)$ since
$S({f_6})=147931/3600-(49/12)\pi^2\simeq0.7911$ and
$S({f_L})=1$.] For the multivariate Student distribution
(\ref{f-Student}), we get $S({f_\nu}) = (1/4) (\nu+m)^2 \{\dot
\Psi(\nu/2) - \dot\Psi[(\nu+m)/2] \}$. The $q$-entropy maximizing
property of the Student distribution can be used to test that the
observed samples are Student distributed, and the estimation of
$S(f)$ then provides information about $\nu$. This finds important
applications, for instance, in financial mathematics; see
\cite{HeydeL2005}.

\section{Main results}
\label{S:main-results}

Let $\rho(x,y)$ denote the Euclidean distance between two points
$x,y$ of $\mathbb{R}^m$ (see Section~\ref{S:further} for an
extension to other metrics). For a given sample $X_1,\ldots,X_N$,
and a given $X_i$ in the sample, from the $N-1$ distances
$\rho(X_i,X_j)$, $j=1,\ldots,N$, $j\neq i$, we form the order
statistics $\rho^{(i)}_{1,N-1} \leq\rho^{(i)}_{2,N-1} \leq\cdots
\leq\rho^{(i)}_{N-1,N-1}$.
Therefore, $\rho^{(i)}_{1,N-1}$ is the nearest-neighbor distance
from the observation $X_i$ to some other $X_j$ in the sample, $j\neq
i$, and similarly, $\rho^{(i)}_{k,N-1}$ is the $k$th
nearest-neighbor distance from $X_i$ to some other $X_j$.

\subsection{R\'enyi and Tsallis entropies}\label{S:RN}

We shall estimate $I_q$, $q\neq1$, by
%
\begin{equation}\label{estimatorIq}
\hat I_{N,k,q} = \frac1N \sum_{i=1}^N (\zeta_{N,i,k})^{1-q} ,
\end{equation}
with
%
\begin{equation}\label{zeta}
\zeta_{N,i,k}=(N-1) C_k V_m \bigl(\rho^{(i)}_{k,N-1}\bigr)^m ,
\end{equation}
where $V_m=\pi^{m/2}/\Gamma(m/2+1)$ is the volume of the unit ball
$\SB(0,1)$ in $\mathbb{R}^m$ and
\[
C_k=\biggl[ \frac{\Gamma(k)}{\Gamma(k+1-q)} \biggr]^{1/(1-q)} .
\]
Note that $I_1=1$ since $f$ is a p.d.f. and that $I_q$ is finite
when $q<0$ only if $f$ is of bounded support. Indeed,
$I_q=\int_{\{x\dvtx f(x)\geq1\}} f^q(x)\,dx+\int_{\{x\dvtx f(x)<1\}} f^q(x)\,dx
>\int_{\{x\dvtx f(x)<1\}} f^q(x)\,dx > \mu_\SL\{x\dvtx f(x)<1\}$, with
$\mu_\SL$ the Lebesgue measure. Also, when $f$ is bounded, $I_q$
tends to the (Lebesgue) measure of its support
$\mu_\SL\{x\dvtx\break f(x)>0\}$ when $q\ra0^+$. Some other properties of
$I_q$ are summarized in Lemma~\ref{L:Iq} of Section~\ref{S:Proofs}.

\begin{rmk} \label{R:Monte-Carlo} When $f$ is known, a Monte Carlo
estimator of $I_q$ based on the
sample $X_1,\ldots,X_N$ is
%
\begin{equation}\label{Iq-MC}
\frac1N \sum_{i=1}^N f^{q-1}(X_i) .
\end{equation}
The nearest-neighbor estimator $\hat I_{N,k,q}$ given by
(\ref{estimatorIq}) could thus also be considered as a plug-in
estimator, $\hat I_{N,k,q} = (1/N) \sum_{i=1}^N [\hat
f_{N,k}(X_i)]^{q-1}$, where $\hat
f_{N,k}(x)=1/\{(N\!-\!1)C_kV_m[\rho_{k+1,N}(x)]^m\}$ with
$\rho_{k+1,N}(x)$ the $(k+1)$th nearest-neighbor distance from $x$ to
the sample. One may notice the resemblance between $\hat f_{N,k}(x)$
and the density function estimator \mbox{${\tilde
f}_{N,k}(x)\!=\!k/\{NV_m[\rho_{k+1,N}(x)]^m\}$} suggested by
Loftsgaarden and Quesenberry \cite{LoftsgaardenQ65}; see also \cite{DevroyeW77,MooreY77}.
\end{rmk}

We suppose that $X_1,\ldots,X_N$, $N\geq2$, are i.i.d.\ with a
probability measure $\mu$ having a density $f$ with respect to the
Lebesgue measure. [However, if $\mu$ has a finite number of singular
components superimposed to the absolutely
continuous component $f$, one can remove all zero distances from the
$\rho_{k,N-1}^{(i)}$ in the computation of
the estimate (\ref{estimatorIq}), which then enjoys the
same properties as in Theorems~\ref{Th:CV of mean}~and~\ref{Th:consistency}, i.e., yields a consistent estimator of the
R\'enyi and Tsallis entropies of the continuous component $f$.] The
main results of the paper are as follows.

\begin{theo}[(Asymptotic unbiasedness)]
\label{Th:CV of mean} The estimator $\hat I_{N,k,q}$ given by
(\ref{estimatorIq}) satisfies
%
\begin{equation}\label{main}
\Ex\hat I_{N,k,q} \ra I_q , \qquad N\ra\infty,
\end{equation}
for $q<1$, provided that $I_q$ given by (\ref{Iq}) exists, and for
any $q\in(1,k+1)$ if $f$ is bounded.
\end{theo}

Under the conditions of Theorem \ref{Th:CV of mean}, $\Ex(1-\hat
I_{N,k,q})/(q-1) \ra H_q$ as $N\ra\infty$, which provides an
asymptotically unbiased estimate of the Tsallis entropy of $f$.

\begin{theo}[(Consistency)]
\label{Th:consistency} The estimator $\hat I_{N,k,q}$ given by
(\ref{estimatorIq}) satisfies
%
\begin{equation}\label{cv}
\hat I_{N,k,q} \raL2 I_q , \qquad N\ra\infty,
\end{equation}
(and thus, $\hat I_{N,k,q} \rap I_q$, $N\ra\infty$) for $q<1$,
provided that $I_{2q-1}$ exists, and for any $q\in(1,(k+1)/2)$
when $k\geq2$ [resp. $q\in(1,3/2)$ when $k=1$] if $f$ is
bounded.
\end{theo}

\begin{coro}\label{C:entropy}
Under the conditions of Theorem \ref{Th:consistency},
%
\begin{equation}\label{hatH}
\hat H_{N,k,q} = (1-\hat I_{N,k,q})/(q-1) \raL2 H_q
\end{equation}
and
%
\begin{equation}\label{hatH*}
\hat H_{N,k,q}^* = \log(\hat I_{N,k,q})/(1-q) \rap H_q^*
\end{equation}
as $N\ra\infty$, which provides consistent estimates of the R\'enyi
and Tsallis entropies of $f$.
\end{coro}

We show the following in the proof of Theorem
\ref{Th:consistency}: when $q<1$ and $I_{2q-1}<\infty$, or
$1<q<(k+2)/2$ and $f$ is bounded,
\[
\Ex(\zeta_{N,i,k}^{1-q}-I_q)^2 \ra\Delta_{k,q}=I_{2q-1}
\frac{\Gamma(k+2-2q)\Gamma(k)}{\Gamma^2(k+1-q)} - I_q^2 , \qquad
N\ra\infty.
\]
Notice that $\lim_{k\ra\infty} \Delta_{k,q}=I_{2q-1}-I_q^2 =
\var[f^{q-1}(X)]=N \var[ \frac1N \times \break \sum_{i=1}^N f^{q-1}(X_i)
]$, that is, the limit of $\Delta_{k,q}$ for $k\ra\infty$
equals $N$ times the variance of the Monte Carlo estimator
(\ref{Iq-MC}) (which forms a lower bound on the variance of an
estimator $I_q$ based on the sample $X_1,\ldots,X_N$).

Under the assumption that $f$ is three times continuously
differentiable $\mu_\SL$-almost everywhere, we can improve Lemma
\ref{L:lebesgue} of Section~\ref{S:Proofs} into
\[
\frac{1}{V_mR^m} \int_{\SB(x,R)} f(z)\,dz = f(x)+ \frac{R^2}{2(m+2)}
\sum_{i=1}^m \frac{\mp^2 f(x)}{\mp x_i^2} + o(R^2) , \qquad R \ra0 ,
\]
which can be used to approximate $F_{N,x,k}(u)-F_{x,k}(u)$ in the
proof of Theorem~\ref{Th:CV of mean}. We thereby obtain an
approximation of the bias $\hat B_{N,k,q} = \Ex\hat I_{N,k,q}-I_q
= \Ex\zeta_{N,1,k}^{1-q}-I_q$, which, after some calculations, can be
written as
\[
\hat B_{N,k,q} = \cases{
\displaystyle\frac{(q-1)(2-q)
I_q}{2N} + \SO(1/N^2),\qquad \mbox{for } m=1, \cr
\displaystyle\frac{q-1}{N} [ (k+1-q) J_{q-2}/(8\pi) + (2-q)I_q/2 ] +
\SO(1/N^{3/2}), \cr
\quad\hspace*{149pt}\mbox{for } m=2, \cr
\displaystyle\frac{q-1}{N^{2/m}} \frac{\Gamma(k+1+2/m-q)}{D_m \Gamma(k+1-q)}
J_{q-1-2/m} + \SO(1/N^{3/m}), \cr
\quad\hspace*{149pt}\mbox{for } m \geq3 ,
}
\]
where $J_\beta=\int_{\mathbb{R}^m} f^\beta(x) (\sum_{i=1}^m
\mp^2 f(x)/\mp x_i^2 )\, dx$ and $D_m=2(m+2)V_m^{2/m}$. For
instance, for $f$ the density of the normal $\SN(0,\ms^2 I_m)$, we
get
\[
J_\beta= - \frac{m}{\ms^2} \frac{1}{(2\pi\ms^2)^{m\beta/2}}
\frac{\beta}{(\beta+1)^{1+m/2}} ,
\]
which is defined for $\beta>-1$. From the expression of the MSE
for $\hat I_{N,k,q}$ given in (\ref{main-consistency}), we
obtain
%
\begin{equation}\label{MSE}
\qquad\quad\Ex(\hat I_{N,k,q}-I_q)^2 = \frac{\Delta_{k,q}}{N} - 2 I_q \hat
B_{N,k,q} \bigl(1+o(1)\bigr) + [
\Ex(\zeta_{N,1,k}^{1-q}\zeta_{N,2,k}^{1-q}) - I_q^2 ] .
\end{equation}
%
Investigating the behavior of the last term requires an
asymptotic approximation for
$F_{N,x,y,k}(u,v)-F_{x,k}(u)F_{y,k}(v)$ (see the proof of
Theorem~\ref{Th:consistency}), which is under current investigation.
Preliminary results for $k=1$ show that the contribution of this
term to the MSE for $\hat I_{N,k,q}$ cannot be ignored in general.
%

\subsection{Shannon entropy} \label{S:S}

For the estimation of $H_1$ ($q=1$), we take the limit of $\hat
H_{N,k,q}$ as $q\ra1$, which gives
%
\begin{equation}\label{hat-H1}
\hat H_{N,k,1} = \frac1N \sum_{i=1}^N \log\xi_{N,i,k},
\end{equation}
with
%
\begin{equation}\label{xi}
\xi_{N,i,k} = (N-1) \exp[-\Psi(k)] V_m
\bigl(\rho^{(i)}_{k,N-1}\bigr)^m,
\end{equation}
where $\Psi(z)=\Gamma'(z)/\Gamma(z)$ is the digamma function
[$\Psi(1) = -\gamma$ with $\gamma\simeq0.5772$ the Euler constant
and, for $k\geq1$ integer, $\Psi(k)=-\gamma+A_{k-1}$ with $A_0=0$
and $A_j=\sum_{i=1}^j 1/i$]; see \cite{KraskovSG2004,Victor2002}
for applications of this estimator in physical sciences. We then
have the following.

\begin{coro}
\label{C:biasH1} Suppose that $f$ is bounded and that $I_{q_1}$
exists for some $q_1<1$. Then $H_1$ exists and the estimator
(\ref{hat-H1}) satisfies $\hat H_{N,k,1} \raL2 H_1$ as $N\ra
\infty$.
\end{coro}


\begin{rmk} \label{R:Sf}
One may notice that $\hat H_{N,k,q}^*$ given by (\ref{hatH*}) is a
smooth function of $q$. Its derivative at $q=1$ can be used as an
estimate of $\dot H_1$ defined by (\ref{dotH1def}).
Straightforward calculations give
\begin{eqnarray*}
\lim_{q\ra1} \frac{d \hat H_{N,k,q}^*}{dq} &=&
\frac{\dot{\Psi}(k)}{2} - \frac{m^2}{2} \frac{1}{N} \sum_{i=1}^N
\Biggl[\log\rho_{k,N-1}^{(i)} - \frac{1}{N} \sum_{j=1}^N \log
\rho_{k,N-1}^{(j)} \Biggr]^2 \\
&=& \frac12 \Biggl[ \dot{\Psi}(k) - \frac{1}{N} \sum_{i=1}^N (\log
\xi_{N,i,k}-\hat H_{N,k,1})^2 \Biggr]
\end{eqnarray*}
and $S(f)=-2 \dot H_1$ can be estimated by
%
\begin{equation}\label{hatSf}
\hat{S}_{{N,k}} = \frac{1}{N} \sum_{i=1}^N (\log
\xi_{N,i,k}-\hat H_{N,k,1})^2 - \dot{\Psi}(k) .
\end{equation}
\end{rmk}

We obtain the following in the proof of Corollary \ref{C:biasH1}:
\[
\Ex( \log\xi_{N,i,k} - H_1)^2 \ra\var[\log f(X)] +
\dot{\Psi}(k) , \qquad N\ra\infty,
\]
with $\dot{\Psi}(z)=d^2 \log\Gamma(z)/dz^2$ [and, for $k$ integer,
$\dot{\Psi}(k)=\sum_{j=k}^\infty1/j^2$]. Note that $\var[\log
f(X)]$ forms a lower bound on the variance of a Monte Carlo
estimation of $H_1$ based on $\log f(X_i)$, $i=1,\ldots,N$, and
that $\dot{\Psi}(k)\ra0$ as $k\ra\infty$.

Similarly to Remark \ref{R:Monte-Carlo}, the estimator $\hat
H_{N,k,1}$ given by (\ref{hat-H1}) could be considered as a
plug-in estimator, $\hat H_{N,k,1} = -(1/N) \sum_{i=1}^N \log
[{\hat f}'_{N,k}(X_i)]$ with ${\hat
f}'_{N,k}(x)=\exp[\Psi(k)]/\{(N-1)V_m[\rho_{k+1,N}(x)]^m\}$. One
may notice that selecting $k$ by likelihood cross-validation based
on the density function estimator suggested by
Loftsgaarden and Quesenberry \cite{LoftsgaardenQ65}, ${\tilde
f}_{N,k}(x)=k/\{NV_m[\rho_{k+1,N}(x)]^m\}$, amounts to maximizing
$-\hat H_{N,k,1}+ \log k - \Psi(k)$, with $\log k -
\Psi(k)=1/(2k)+1/(12 k^2)+\SO(1/k^4)$, $k\ra\infty$. In our
simulations this method always tended to select $k=1$; replacing
${\tilde f}_{N,k}(x)$ by ${\hat f}'_{N,k}(x)$, or by ${\hat
f}_{N,k}(x)$ of Remark \ref{R:Monte-Carlo}, does not seem to yield
a valid selection procedure for $k$ either.

Let $\tilde H_{N,k,1}$ be the plug-in estimator of $H_1$ based on
${\tilde f}_{N,k}$ defined by $\tilde H_{N,k,1} = -(1/N)
\sum_{i=1}^N \log[{\tilde f}_{N,k}(X_i)]$. Then, under the
conditions of Corollary \ref{C:biasH1}, we obtain that
$\lim_{N\ra\infty} \Ex\tilde H_{N,k,1} = H_1 + \Psi(k)-\log k$
(since $\tilde H_{N,k,1}=\hat H_{N,k,1} + \Psi(k)-\log k +
\log[N/(N-1)]$). Under the additional assumption on $f$ that it
belongs to the class $\SF$ of uniformly continuous p.d.f.\
satisfying $0<c_1\leq f(x) \leq c_2< \infty$ for some constants
$c_1, c_2$, we obtain the uniform and almost sure convergence of
$\hat H_{N,k,1}$ to $H_1(f)$ over the class $\SF$, provided that
$k=k_N \ra\infty$, $k_N/N\ra0$ and $k_N/\log N \ra\infty$ as
$N\ra\infty$; see the results of Devroye and Wagner
\cite{DevroyeW77} on the strong uniform consistency of ${\tilde
f}_{N,k}$. Notice that the choice of $k$ proposed by Hall, Park and
Samworth \cite{HallPS2004} for nearest-neighbor classification does
not satisfy these conditions.

\subsection{Relative entropy and divergences}\label{S:divergences}

In some situations the statistical distance between distributions
can be estimated through the computation of entropies, so that the
method of $k$th nearest-neighbor distances presented above can be
applied straightforwardly. For instance, the $q$-Jensen difference
\[
J_q^\beta(f,g) = H_q^*[\beta f+ (1-\beta)g] - [\beta
H_q^*(f)+(1-\beta)H_q^*(g)] , \qquad 0 \leq\beta\leq1 ,
\]
(see, e.g., \cite{Basseville96}) can be estimated if we have three
samples, respectively distributed according to $f$, $g$ and $\beta
f+ (1-\beta)g$. Suppose that we have one sample $S_i$
($i=1,\ldots,s$) of i.i.d.\ variables generated from $f$ and one
sample $T_j$ ($j=1,\ldots,t$) of i.i.d.\ variables generated from
$g$ with $s$ and $t$ increasing at a constant rate as a function
of $N=s+t$. Then, $H_q^*(f)$ and $H_q^*(g)$ can be estimated
consistently from the two samples when $N\ra\infty$; see Corollary
\ref{C:entropy}. Also, as $N\ra\infty$, the estimator $\hat
H_{N,k,q}^*$ based on the sample $X_1,\ldots,X_N$ with $X_i=S_i$
($i=1,\ldots,s$) and $X_i=T_{i-s}$ ($i=s+1,\ldots,N$) converges to
$H_q^*[\beta f+ (1-\beta)g]$, with $\beta=s/N$, and $J_q^\beta$
can therefore be estimated consistently from the two samples. This
situation is encountered, for instance, in the image matching
problem presented in \cite{NeemuchwalaHC2005}, where entropy is
estimated through the random graph approach of Redmond and Yukich \cite{RedmondY96}.
As shown below, some other types of distances or divergences, that
are not expressed directly through entropies, can also be
estimated by the nearest-neighbor method.

Let $K(f,g)$ denote the Kullback--Leibler relative entropy,
%
\begin{equation}\label{Kullback-L}
K(f,g) = \int_{\mathbb{R}^m} f(x) \log\frac{f(x)}{g(x)}
\,dx = \breve{H}_1 - H_1 ,
\end{equation}
where $H_1$ is given by (\ref{Shannon}) and $\breve{H}_1 = -
\int_{\mathbb{R}^m} f(x) \log g(x) \, dx$. Given $N$ independent
observations $X_1,\ldots,X_N$ distributed with the density $f$ and
$M$ observations $Y_1,\ldots,Y_M$ distributed with $g$, we wish to
estimate $K(f,g)$. The second term $H_1$ can be estimated by
(\ref{hat-H1}), with asymptotic properties given by Corollary
\ref{C:biasH1}. The first term $\breve{H}_1$ can be estimated in a
similar manner, as follows: given $X_i$ in the sample,
$i\in\{1,\ldots,N\}$, consider $\breve{\rho}(X_i,Y_j)$,
$j=1,\ldots,M$, and the order statistics $\breve{\rho}^{(i)}_{1,M}
\leq\breve{\rho}^{(i)}_{2,M} \leq\cdots\leq
\breve{\rho}^{(i)}_{M,M}$, so that $\breve{\rho}^{(i)}_{k,M}$ is the
$k$th nearest-neighbor distance from $X_i$ to some $Y_j$,
$j\in\{1,\ldots,M\}$. Then, one can prove, similarly to Corollary
\ref{C:biasH1}, that
%
\begin{equation}\label{breveHNMk}
\breve{H}_{N,M,k}= \frac{1}{N} \sum_{i=1}^N \log\bigl\{ M
\exp[-\Psi(k)] V_m \bigl(\breve{\rho}^{(i)}_{k,M}\bigr)^m\bigr\}
\end{equation}
is an asymptotically unbiased and consistent estimator of
$\breve{H}_1$ (when now both $N$ and $M$ tend to infinity) when $g$
is bounded and
%
\begin{equation}\label{Jq}
J_q = \int_{\mathbb{R}^m} f(x) g^{q-1}(x) \,dx
\end{equation}
exists for some $q<1$. The difference
%
\begin{eqnarray} \label{Estimate_KL}
\breve{H}_{N,M,k} - \hat H_{N,k,1} &=&
m \log\Biggl[\prod_{i=1}^N \breve{\rho}^{(i)}_{k,M}\Biggr]^{1/N} +
\log M - \Psi(k) + \log V_m \nonumber\\
&&{} - m \log\Biggl[\prod_{i=1}^N \rho^{(i)}_{k,N}\Biggr]^{1/N} - \log
(N-1) + \Psi(k) - \log V_m \\
&=& m \log\Biggl[\prod_{i=1}^N
\frac{\breve{\rho}^{(i)}_{k,M}}{\rho^{(i)}_{k,N}}\Biggr]^{1/N}+
\log\frac{M}{N-1}\nonumber
\end{eqnarray}
thus gives an asymptotically unbiased and consistent estimator of
$K(f,g)$. Obviously a similar technique can be used to estimate the
(symmetric) Kullback--Leibler divergence $K(f,g)+K(g,f)$. Note, in
particular, that when $f$ is unknown and only the sample
$X_1,\ldots,X_N$ is available while $g$ is known, then the term
$\breve{H}_1$ in $K(f,g)$ can be estimated either by
(\ref{breveHNMk}) with a sample $Y_1,\ldots,Y_M$ generated from $g$,
with $M$ taken arbitrarily large, or more simply by the Monte Carlo
estimator
%
\begin{equation}\label{breveH1N}
\breve{H}_{1,N}(g) = -\frac{1}{N} \sum_{i=1}^N \log g(X_i) ,
\end{equation}
the term $H_1$ being still estimated by (\ref{hat-H1}). This forms
an alternative to the method by Broniatowski \cite{Broniatowski2003}. Compared
to the method by Jim{\'e}nez and Yukich \cite{JimenezY2002} based on Voronoi
tessellations (see also \cite{Miller2003} for a Voronoi-based
method for Shannon entropy estimation), it does not require any
computation of multidimensional integrals. In some applications
one wishes to optimize $K(f,g)$ with respect to $g$ that belongs
to some class $G$ (possibly parametric), with $f$ fixed. Note that
only the first term $\breve{H}_1$ of (\ref{Kullback-L}) then needs
to be estimated. [Maximum likelihood estimation, with $g=g_\theta$
in a parametric class, is a most typical example: $\theta$ is then
estimated by minimizing $\breve{H}_{1,N}(g_\theta)$; see
(\ref{breveH1N}).]

The Kullback--Leibler relative entropy can be used to construct a
measure of mutual information (MI) between statistical
distributions (see \cite{KraskovSG2004}) with applications in
image \cite{NeemuchwalaHC2005,ViolaW95} and signal processing
\cite{LearnedMillerF2003}. Let $a_i$ and $b_i$ denote the gray
levels of pixel $i$ in two images $A$ and $B$ respectively,
$i=1,\ldots,N$. The image matching problem consists in finding an
image $B$ in a data base that resembles a given reference image
$A$. The MI method corresponds to maximizing $K(f,f_xf_y)$, with
$f$ the joint density of the pairs $(a_i,b_i)$ and $f_x$
(resp. $f_y$) the density of gray levels in image $A$
(resp. $B$). We have
$K(f,f_xf_y)=H_1(f_x)+H_1(f_y)-H_1(f)$, where each term can be
estimated by (\ref{hat-H1}) from one of the three samples $(a_i)$,
$(b_i)$ or $(a_i,b_i)$ (but $A$ being fixed, only the last two
terms need be estimated).

Another example of statistical distance between distributions is
given by the following nonsymmetric Bregman distance
%
\begin{eqnarray}\label{Bregman}
D_q(f,g) = \int_{\mathbb{R}^m} \biggl[ g^q(x)+\frac{1}{q-1} f^q(x) -
\frac{q}{q-1} f(x)g^{q-1}(x)\biggr]\, dx , \nonumber\\[-8pt]\\[-8pt]
\eqntext{q\neq1 ,}
\end{eqnarray}
or its symmetrized version
\begin{eqnarray*}
K_q(f,g)&=&\frac{1}{q} [D_q(f,g)+D_q(g,f)]\\
& = & \frac{1}{q-1}
\int_{\mathbb{R}^m} [f(x)-g(x)] [f^{q-1}(x)-g^{q-1}(x)] \,dx;
\end{eqnarray*}
see, for example, \cite{Basseville96}. Given $N$ independent observations
from $f$ and $M$ from $g$, the first and second terms in
(\ref{Bregman}) can be estimated by using (\ref{estimatorIq}). In
the last term, the integral $J_q$ given by (\ref{Jq}) can be
estimated by $\hat I_{N,M,k,q} = (1/N) \sum_{i=1}^N \{M C_k V_m
(\breve{\rho}^{(i)}_{k,M})^m \}^{1-q}$. Similarly to Theorem~\ref{Th:CV of mean},
 $\hat I_{N,M,k,q}$ is asymptotically unbiased,
$N,M \ra\infty$, for $q<1$ if $J_q$ exists and for any
$q\in(1,k+1)$ if $g$ is bounded. We also obtain a property similar
to Theorem \ref{Th:consistency}: $\hat I_{N,M,k,q}$ is a consistent
estimator of $J_q$, $N,M\ra\infty$, for $q<1$ if $J_{2q-1}$ exists
and for any $q\in(1,(k+2)/2)$ if $g$ is bounded. (Notice, however,
the difference with Theorem~\ref{Th:consistency}: when $q>1$ the
cases $k=1$ and $k\geq2$ need not be distinguished for the
estimation of $J_q$ and the upper bound on the admissible values for
$q$ is slightly larger than in Theorem \ref{Th:consistency}.)

%
\begin{figure}[b]

\includegraphics{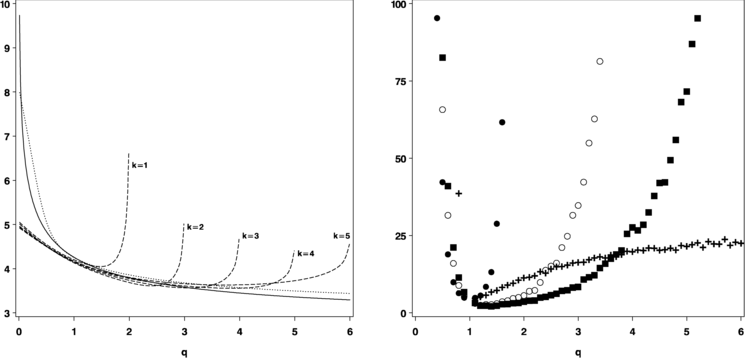}

\caption{Behavior of estimators of entropy for samples from the
normal distribution $\SN(0,I_3)$ in $\mathbb{R}^3$ ($N=1 000$).
[Left] $H_q^*$ (solid line), $\hat H_{N,k,q}^*$ (dashed lines) and
$\tilde H^*_{N,q}$ obtained through a kernel estimation of $f$
(dotted line) as functions of $q$. [Right] $N=1 000$ times the
empirical MSE for $\hat H_{N,k,q}$ [$k=1$ (dots), $k=3$ (circles),
$k=5$ (squares)] and for $\tilde H^*_{N,q}$ (plus) as a function
of $q$ and computed over 1 000 independent samples. }
\label{F:Renyi_normal_m3_q}
\end{figure}

\newpage

\section{Examples} \label{S:examples}
\subsection{Influence of $k$}
Figure~\ref{F:Renyi_normal_m3_q} (left) presents $H_q^*$ as a
function of $q$ (solid line) for the normal distribution
$\SN(0,I_3)$ in $\mathbb{R}^3$, together with estimates $\hat
H_{N,k,q}^*$ for $k=1,\ldots,5$ obtained from a single sample of
size $N=1 000$. Note that $\hat H_{N,k,q}^*$ is defined only for
$q<k+1$ and quickly deviates from the theoretical value $H_q^*$
when $q>(k+1)/2$ or $q<1$ (the difficulties for $q$ small being
due to $f$ having unbounded support). For comparison, we also
compute a plug-in estimate of $H_q^*$ obtained through a
(cross-validated) kernel density estimate of $f$. Define $\tilde
H^*_{N,q}=\log(\tilde I_{N,q})/(1-q)$ and $\tilde
I_{N,q}=(1/N)\sum_{i=1}^N \tilde f_{N,i}^{q-1}(X_i)$ with $\tilde
f_{N,i}(x)=[(N-1)h^m(2\pi)^{m/2}]^{-1} \sum_{l=1, l\neq i}^N
\exp\{-\|x-X_l\|^2/(2h^2)\}$, a $m$-variate cross-validated
kernel estimator of $f$. No special care is taken for the choice
of $h$ and we simply use the value that minimizes the asymptotic
mean integrated squared error for the estimation of $f$, that is,
$h=[4/(m+2)]^{1/(m+4)} N^{-1/(m+4)}$ with $m=3$; see \cite{Scott92},
page 152. The evolution of $\tilde H^*_{N,q}$ as a function
of $q$ is plotted in dotted-line on Figure~\ref{F:Renyi_normal_m3_q}
(left): although the situation is favorable to kernel density
estimation, $k$th nearest neighbors give a better estimation of
$H_q^*$ for $q>1$ and $k$ large enough.
Figure \ref{F:Renyi_normal_m3_q} (right) shows $N$ times the
empirical mean-squared error (MSE) $\Ex(\hat H^*_{N,k,q} -
H^*_{N,q})^2$ ($k=1,3,5$) as a function of $q$ using 1 000
independent repetitions. The results for $N$ times the MSE
$\Ex(\tilde H^*_{N,q} - H^*_{N,q})^2$ for the plug-in estimator
are also shown. The figure indicates that the $k$th nearest
neighbor estimator with $k$ satisfying $q<(k+1)/2$ is favorable in
comparison to the plug-in estimator (for $q>1$ values of $k$
larger than 1 are preferable, whereas $k=1$ is preferable, for
$q<1$).

\begin{figure}[b]

\includegraphics{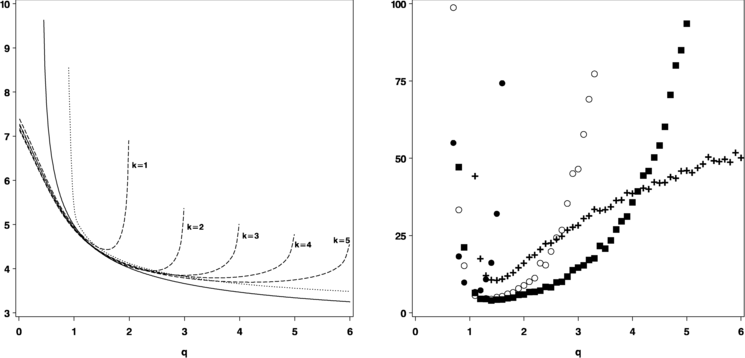}

\caption{Same information as in Figure~\protect\ref{F:Renyi_normal_m3_q} but
for the Student distribution $T(\nu,\Sigma,\mu)$ in $\mathbb{R}^3$
with 4 degrees of freedom ($\Sigma=I_3$, $\mu=0$, $N=1 000$).
} \label{F:Renyi_student_m3_q}
\end{figure}

Similar results hold for the Student distribution for
$T(\nu,\Sigma,\mu)$ in $\mathbb{R}^3$ with 4 degrees of freedom,
$\Sigma=I_3$ and $\mu=0$; see Figure~\ref{F:Renyi_student_m3_q}. In
selecting $k$ for $\hat H_{N,k,q}^*$, large values of $k$ are
still generally preferable when $q>1$.

%

\begin{figure}

\includegraphics{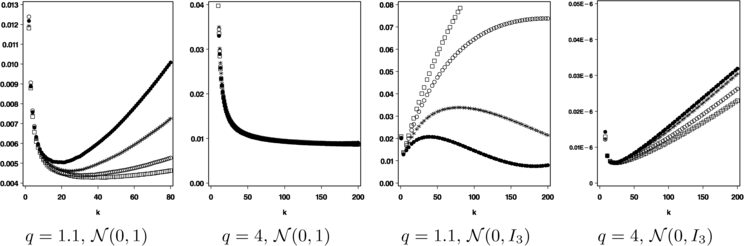}

\caption{$N$ times the empirical MSE for $\hat I_{N,k,q}$ as a
function of $k$ (10 000 independent repetitions), for $f$ the
density of the standard normal $\SN(0,1)$ and $\SN(0,I_3)$ in
$\mathbb{R}^3$ for varying $N$ \{$N=1 000$ (dots), 2 000
(stars), 5 000 (circles) and 10 000 (squares)\} and $q=1.1$ and
$q=4$.} \label{F:Iq_N}
\end{figure}

\begin{figure}[b]

\includegraphics{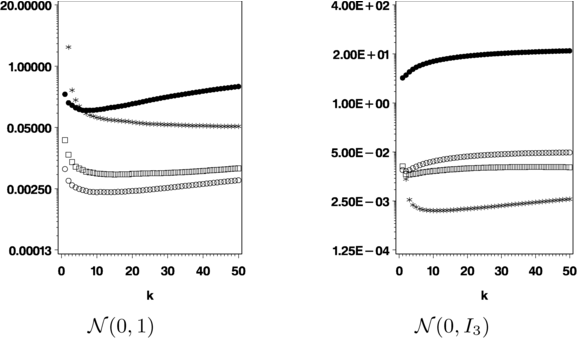}

\caption{$N$ times the empirical MSE for $\hat I_{N,k,q}$ as a
function of $k$ (10 000 independent repetitions), for $f$ the
density of the standard normal $\SN(0,1)$ and $\SN(0,I_3)$ in
$\mathbb{R}^3$ for varying $q$ \{$q=0.75$ (dots), $q=0.95$
(circles), $q=1.1$ (squares) and $q=2$ (stars)\} and
$N=1 000$.} \label{F:Iq_q}
\end{figure}

At this stage, the optimal selection of $k$ in $\hat I_{N,k,q}$
depending on $q$ and $N$ remains an open issue (see
Sections~\ref{S:S} and \ref{S:further}). We repeated a series of
intensive simulations to see how the MSE $\Ex(\hat I_{N,k,q} -
I_q)^2$ evolves when $k$ varies, for different choices of $N$, $q$
and $m$.
%
Figure~\ref{F:Iq_N} shows the influence of $N$ on the MSE for $\hat
I_{N,k,q}$ for different values of $q$ using 10 000 independent
repetitions, for $f$ the density of the standard normal $\SN(0,1)$
and the normal $\SN(0,I_3)$. For both $m=1$ and $m=3$ changes in
$N$ appear to have a greater influence on $N$ times the MSE for
$q=1.1$ in comparison to $q=4$. In particular, the figure
indicates that for $m=3$ and $q=1.1$ the MSE decreases more slowly
than $1/N$. Figure~\ref{F:Iq_q} shows the influence of $q$ on $N$
times the MSE for $\hat I_{N,k,q}$ as $k$ varies.

Although our simulations do not reveal a precise rule for choosing
$k$, they indicate that this choice is not critical for practical
applications: taking $k$ between~5~and~10 for $q\leq2$ and
increasing from 10 to 20 for $q$ from 2 to 4 gives reasonably good
results for the cases we considered.

%


%

%
%


\subsection{Information spectrum, estimation of $\var[\log f(X)]$}
\label{S:Ex-IS}

We use the method suggested in Remark \ref{R:Sf} and estimate
$S(f)=\var[\log f(X)]$ by $\hat S_{N,1}$ given by (\ref{hatSf})
from a sample of 50 000 data generated with the Student
distribution with 5 degrees of freedom. $S(f_\nu)$ is a decreasing
function of $\nu$ and $S({f_4})\simeq0.9661$, $S({f_5})\simeq
0.8588$, $S({f_6})\simeq0.7911$; see Section~\ref{S:spectrum}. The
empirical mean and standard deviation of $\hat{S}_{{N,1}}$
obtained from 10 000 independent repetitions are 0.8578 and
0.0269 respectively, indicating that $\nu$ can be correctly
estimated in this way.

%

\subsection{Estimation of Kullback--Leibler divergence}

We use the same Student data as in \ref{S:Ex-IS} and estimate the
Kullback--Leibler relative entropy $K(f,f_\nu)$ given by
(\ref{Kullback-L}), using (\ref{breveH1N}) for the estimation of
$\breve{H}_1$ and (\ref{hat-H1}) for the estimation of $H_1$, the
entropy of $f$.
The empirical means of the
divergences estimated for $\nu=1,\ldots,8$ in 10 000 independent
repetitions are $0.1657$, $0.0440$, $0.0119$, $0.0021$, $0.0000$,
$0.0012$, $0.0038$ and $0.0069$ [the empirical standard deviations
are rather large, approximately 0.0067 for each $\nu$, but the
minimum is at $\nu=5$ in all the 10 000 cases
---notice that the dependence in $\nu$ is only through the term
(\ref{breveH1N}) where $f_\nu$ is substituted for $g$]. Again,
$\nu$ is correctly estimated in this way.


%

\subsection{$q$-entropy maximizing distributions}
\begin{figure}

\includegraphics{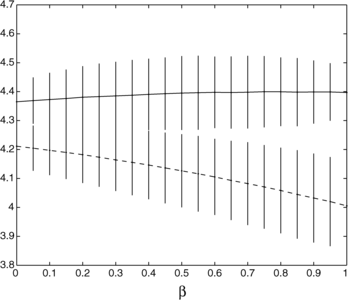}

\caption{Empirical means of $\hat H_{N,3,0.75}^*$ (solid line) and
$\hat H_{N,3,1}$ (dashed line) and two standard deviations
(vertical bars) in a mixture of Student and normal distributions
as functions of the mixture coefficient $\beta$ for $N=500$
(\textup{1 000} independent repetitions).} \label{F:mixture_m3}
\end{figure}

We generate $N=500$ i.i.d.\ samples from the mixture of the
three-dimensional Student distribution $T(\nu,(\nu-2)/\nu I_3,0)$
with $\nu=5$ and the normal distribution $\SN(0,I_3)$, with
relative weights $\beta$ and $1-\beta$. The covariance matrix of
both distributions is the identity $I_3$, the Student distribution
is $q$-entropy maximizing for $q=1-2/(\nu+m)=0.75$ (see
Section~\ref{S:entropymaximizing}) and the normal distribution maximizes
Shannon entropy ($q=1$). Figure~\ref{F:mixture_m3} presents a plot
of $\hat H_{N,k,q}^*$ for $q=0.75$ and $\hat H_{N,k,1}$ as
functions of the mixture coefficient $\beta$; both use $k=3$ and
are averaged over 1 000 repetitions, the vertical bars indicate
two empirical standard deviations. [The values of $H_{0.75}^*$
estimated by plug-in using the kernel estimator $\tilde
f_{N,i}(x)$ of Example 1 are totally out of the range for Student
distributed variables due to the use of a nonadaptive bandwidth.]

%

\section{Related results and further developments}
\label{S:further}
The paper by Jim{\'e}nez and Yukich \cite{JimenezY2002} gives a method for estimating
statistical distances between distributions with densities $f$ and
$g$ based on Voronoi tessellations. Given an i.i.d.\ sample from
$f$, it relies on the comparison between the Lebesgue measure
(volume) and the measure for $g$ of the Voronoi cells (polyhedra)
constructed from the sample. Voronoi tessellations are also used in
\cite{Miller2003} to estimate the Shannon entropy of $f$ based on
an i.i.d.\ sample. The method requires the computation of the
volumes of the Voronoi cells and no asymptotic result is given.
Comparatively, the method based on nearest neighbors does not
require any computation of (multidimensional) integrals. A possible
motivation for using Voronoi tessellations could be the natural
adaptation to the shape of the distribution. One may then notice
that the metric used to compute nearest-neighbor distances can be
adapted to the observed sample: for $X_1,\ldots,X_N$, a sample having
a nonspherical distribution, its empirical covariance matrix
$\hat\Sigma_N$ can be used to define a new metric through
$\|x\|_{\hat\Sigma_N}^2=x\TT\hat\Sigma_N^{-1} x$, the volume $V_m$
of the unit ball in this metric becoming
$|\hat\Sigma_N|^{1/2}\pi^{m/2}/\Gamma(m/2+1)$.

$\sqrt{N}$-consistency of an estimator of $H_1$ based on
nearest-neighbor distances ($k=1$) is proved by Tsybakov and van~der
Meulen \cite{TsybakovvdM96} for $m=1$ and sufficiently regular
densities $f$ with unbounded support using a truncation argument. On
the other hand, $\sqrt{N}$-consistency of the estimator $\hat
I_{N,k,q}$ is still an open issue (notice that the bias approximations
of Section~\ref{S:RN} indicate that it does not hold for large $m$). As
for the case of spacing methods, where the spacing can be taken as an
increasing function of the sample size $N$ (see, e.g.,
\cite{HAll86,vanEs92,Vasicek76}) it might be of interest to let $k=k_N$
increase with $N$; see also \cite{Song2000} and Section~\ref{S:S}.
Properties of nearest-neighbor distances with $k_N\ra\infty$ are
considered, for instance, by Devroye and Wagner \cite{DevroyeW77},
Liero \cite{Liero93}, Loftsgaarden and Quesenberry
\cite{LoftsgaardenQ65} and Moore and Yackel \cite{MooreY77}. The
derivation of an estimate of the asymptotic mean-squared error of the
estimator could be used in a standard way to construct a rule for
choosing $k$ as a function of $q$, $m$ and $N$ (see Sections~\ref{S:RN}
and \ref{S:S}). Numerical simulations indicate, however, that this
choice is not as critical as that of the bandwidth in a kernel density
estimator used for plug-in entropy estimation; see
Section~\ref{S:examples}.\looseness=1

A central limit theorem for functions $h(\rho)$ of nearest-neighbor
distances is obtained by Bickel and Breiman \cite{BickelB83} for $k=1$ and by
Penrose \cite{Penrose2000} for $k=k_N\ra\infty$ as $N\ra\infty$. However,
their results do not apply to unbounded functions of $\rho$, such as
$h(\rho)=\rho^{m(1-q)}$ [see (\ref{estimatorIq})], or
$h(\rho)=\log(\rho)$ [see (\ref{hat-H1})]. Conditions for the
asymptotic normality of $\hat I_{N,k,q}$ are under current
investigation.\looseness=1

\section{Proofs}\label{S:Proofs}

The following lemma summarizes some properties of $I_q$.

\begin{lemma}\label{L:Iq}

\begin{longlist}[(iii)]
\item[(i)] If $f$ is bounded, then $I_q<\infty$ for any $q>1$.

\item[(ii)] If $I_q<\infty$ for some $q<1$, then
$I_{q'}<\infty$ for any $q'\in(q,1)$.

\item[(iii)] If $f$ is of finite support, $I_q<\infty$ for any
$q\in[0,1)$.
\end{longlist}
\end{lemma}

\begin{pf}
\begin{longlist}
\item[(i)] If $f(x)<\bar f$ and $q>1$, $I_q=\int_{f\leq1} f^q + \int_{f>1}
f^q \leq\int_{f\leq1} f + \bar f^q \int_{f>1} f < \infty$.

\item[(ii)] If $q<q'<1$, $I_{q'} = \int_{f\leq1} f^{q'} + \int_{f>1}
f^{q'} \leq\int_{f\leq1} f^q + \int_{f>1} f < \infty$ if
$I_q<\infty$.

\item[(iii)] If $\mu_\SS=\mu_\SL\{x\dvtx f(x)>0\}<\infty$ and $0\leq q<1$,
$I_q=\int_{f\leq1} f^q + \int_{f>1} f^q \leq\mu_\SS+ \int_{f>1} f
< \infty$.\hspace*{2pt}\qed %
\end{longlist} \noqed
\end{pf}

The proofs of Theorems \ref{Th:CV of mean} and
\ref{Th:consistency} use the following lemmas.

\begin{lemma}[{[}Lebesgue (1910){]}] \label{L:lebesgue} If $g\in
L_1(\mathbb{R}^m)$, then for any sequence of open balls $\SB(x,R_k)$
of radius tending to zero as $k\ra\infty$ and for $\mu_\SL$-almost
any $x\in\mathbb{R}^m$,
\[
\lim_{k\ra\infty} \frac{1}{V_m R_k^m} \int_{\SB(x,R_k)} g(t)\,dt =
g(x) .
\]
\end{lemma}

\begin{lemma}
\label{L:moment-F} For any $\beta>0$,
%
\begin{equation}
\int_0^\infty x^\beta F(dx) = \beta\int_0^\infty
x^{\beta-1} [1-F(x)] \,dx \label{moment-F-Feller}
\end{equation}
and
%
\begin{equation}
\int_0^\infty x^{-\beta} F(dx) = \beta\int_0^\infty
x^{-\beta-1} F(x)\, dx \label{moment-F-new} ,
\end{equation}
in the sense that if one side converges so does the other.
\end{lemma}

\begin{pf}
See \cite{Feller66}, volume 2, page 150, for
(\ref{moment-F-Feller}). The proof is similar for
(\ref{moment-F-new}). Define $\ma=-\beta<0$ and $I_{a,b}=\int_a^b
x^\ma F(dx)$ for some $a, b$, with $0<a<b<\infty$. Integration
by parts gives $I_{a,b}=[b^\ma F(b)-a^\ma F(a)] - \ma\int_a^b
x^{\ma-1} F(x) \,dx$ and, since $\ma<0$, $\lim_{b\ra\infty}
I_{a,b}=I_{a,\infty} = -a^\ma F(a) - \ma\int_a^\infty x^{\ma-1}
F(x) \,dx < \infty$. Suppose that $\int_0^\infty x^{-\beta}
F(dx)=J<\infty$. It implies $\lim_{a\ra0^+} I_{0,a}=0$ and, since
$I_{0,a}>a^\ma F(a)$, $\lim_{a\ra0^+} a^\ma F(a)=0$. Therefore,
$\lim_{a\ra0^+} - \ma\int_a^\infty x^{\ma-1} F(x)\, dx = J$.

Conversely, suppose that $\lim_{a\ra0^+} - \ma\int_a^\infty
x^{\ma-1} F(x) dx = J < \infty$. Since\break $I_{a,\infty}<- \ma
\int_a^\infty x^{\ma-1} F(x) \,dx$, $\lim_{a\ra0^+} I_{a,\infty}
= J$.
\end{pf}

%

\subsection{\texorpdfstring{Proof of Theorem \protect\ref{Th:CV of mean}}{Proof of Theorem 3.1}}

Since the $X_i$'s are i.i.d.,
\[
\Ex\hat I_{N,k,q}=\Ex\zeta_{N,i,k}^
{1-q}=\Ex[\Ex(\zeta_{N,i,k}^{1-q}|X_i=x)] ,
\]
where the random variable $\zeta_{N,i,k}$ is defined by
(\ref{zeta}). Its distribution function conditional to $X_i=x$ is
given by
\[
F_{N,x,k}(u) = \Pr(\zeta_{N,i,k} < u | X_i=x)
= \Pr\bigl[ \rho^{(i)}_{k,N-1} < R_{N}(u) | X_i=x \bigr],
\]
where
%
\begin{equation}\label{RN}
R_{N}(u)=\{u/[(N-1)V_m C_k]\}^{1/m} .
\end{equation}
Let $\SB(x,r)$ be the open ball of center $x$ and radius $r$. We
have
\begin{eqnarray*}
F_{N,x,k}(u) &=& \Pr\{k \mbox{ elements or more } \in\SB
[x,R_{N}(u)] \} \\
&=& \sum_{j=k}^{N-1}
\pmatrix{
N-1 \cr j}
p_{N,u}^j (1-p_{N,u})^{N-1-j} \\
&=& 1-\sum_{j=0}^{k-1}
\pmatrix{ N-1 \cr j
}
p_{N,u}^j (1-p_{N,u})^{N-1-j},
\end{eqnarray*}
where $p_{N,u} = \int_{\SB[x,R_N(u)]} f(t) \,dt$. From the Poisson
approximation of binomial distribution, Lemma \ref{L:lebesgue}
gives
\[
F_{N,x,k}(u) \ra F_{x,k}(u)=1-\exp(-\ml u) \sum_{j=0}^{k-1}
\frac{(\ml u)^j}{j!}
\]
when $N\ra\infty$ for $\mu$-almost any $x$, with $\ml=f(x)/C_k$,
that is, $F_{N,x,k}$ tends to the Erlang distribution $F_{x,k}$,
with p.d.f.\ $f_{x,k}(u)= [\ml^k u^{k-1} \exp(-\ml u)]/\Gamma(k)$.
Direct calculation gives
\[
\int_0^\infty u^{1-q} f_{x,k}(u) \,du = \frac{\Gamma(k+1-q)}{\ml^{1-q}
\Gamma(k)} = f^{q-1}(x)
\]
for any $q<k+1$.

Suppose first that $q<1$ and consider the random variables $(U,X)$
with joint p.d.f.\ $f_{N,x,k}(u)f(x)$ on
$\mathbb{R}\times\mathbb{R}^m$, where
$f_{N,x,k}(u)=dF_{N,x,k}(u)/du$. The function $u \ra u^{1-q}$ is
bounded on every bounded interval and the generalized Helly--Bray
Lemma (see \cite{Loeve77}, page 187) implies
\begin{eqnarray*}
\Ex\hat I_{N,k,q} &=& \int_{\mathbb{R}^m} \int_0^\infty u^{1-q}
f_{N,x,k}(u) f(x) \,du \,dx\\
& \ra&\int_{\mathbb{R}^m} f^q(x) \,dx
= I_q , \qquad N \ra\infty,
\end{eqnarray*}
which completes the proof.

Suppose now that $1<q<k+1$. Note that from Lemma \ref{L:Iq}(i)
$I_q<\infty$. Consider
\[
J_N=\int_0^\infty u^{(1-q)(1+\delta)} F_{N,x,k}(du) .
\]
We show that $\sup_N J_N <\infty$ for some $\delta>0$. From Theorem
2.5.1 of  Bierens \cite{Bierens94}, page 34, it implies
%
\begin{eqnarray}\label{main2}
\hspace*{18pt}z_{N,k}(x) = \int_0^\infty u^{1-q} F_{N,x,k}(du) \ra
z_k(x)=\int_0^\infty u^{1-q} F_{x,k}(du) = f^{q-1}(x) ,\nonumber\\[-8pt]\\[-8pt]
\eqntext{N\ra\infty}
\end{eqnarray}
for $\mu$-almost any $x$ in $\mathbb{R}^m$.


Define $\beta=(1-q)(1+\delta)$, so that $\beta<0$, and take
$\delta<(k+1-q)/(q-1)$ so that $\beta+k>0$. From
(\ref{moment-F-new}),
%
\begin{eqnarray}
\label{JN-new}
J_N &=& -\beta\int_0^\infty u^{\beta-1} F_{N,x,k}(u) \,du \nonumber
\\
&=& -\beta\int_0^1 u^{\beta-1} F_{N,x,k}(u) \,du   -\beta\int
_1^\infty u^{\beta-1} F_{N,x,k}(u)
\,du \nonumber\\[-8pt]\\[-8pt]
& \leq&-\beta\int_0^1 u^{\beta-1} F_{N,x,k}(u) \,du   -\beta\int
_1^\infty u^{\beta-1}
\,du \nonumber\\
& =& 1 -\beta\int_0^1 u^{\beta-1} F_{N,x,k}(u) \,du .\nonumber
\end{eqnarray}
Since $f(x)$ is bounded, say, by $\bar f$, we have $\forall
x\in\mathbb{R}^m , \ \forall u\in\mathbb{R} , \ \forall N , \
p_{N,u} \leq\bar f V_m [R_N(u)]^m = \bar f u/[(N-1)C_k]$. It
implies
\begin{eqnarray*}
\frac{F_{N,x,k}(u)}{u^k} &\leq& \sum_{j=k}^{N-1}
\pmatrix{ N-1 \cr j}
\frac{\bar f^j u^{j-k}}{C_k^j (N-1)^j} \\
&\leq& \sum_{j=k}^{N-1} \frac{\bar f^j u^{j-k}}{C_k^j j!}
= \frac{\bar f^k}{C_k^k k!} + \sum_{j=k+1}^{N-1} \frac{\bar f^j
u^{j-k}}{C_k^j j!}\\
& \leq & \frac{\bar f^k}{C_k^k k!} + \frac{\bar f^k}{C_k^k} \sum
_{i=1}^{N-k-1} \frac{\bar f^i u^i}{C_k^i i!} \\
&\leq & \frac{\bar f^k}{C_k^k k!} + \frac{\bar f^k}{C_k^k}
\sum_{i=1}^\infty\frac{\bar f^i u^i}{C_k^i i!} = \frac{\bar
f^k}{C_k^k k!} + \frac{\bar f^k}{C_k^k} \biggl\{\exp\biggl[
\frac{\bar f u}{C_k} \biggr]-1\biggr\} ,
\end{eqnarray*}
and thus, for $u<1$,
%
\begin{equation}\label{Uk}
\frac{F_{N,x,k}(u)}{u^k} < U_k = \frac{\bar f^k}{C_k^k k!} +
\frac{\bar f^k}{C_k^k} \biggl\{\exp\biggl[ \frac{\bar f}{C_k}
\biggr]-1\biggr\} .
\end{equation}
Therefore, from (\ref{JN-new}),
%
\begin{equation}\label{resT1}
J_N \leq1 - \beta U_k \int_0^1 u^{k+\beta-1} \,du = 1 -
\frac{\beta U_k}{k+\beta} < \infty,
\end{equation}
which implies (\ref{main2}). Now we only need to prove that
\[
\int_{\mathbb{R}^m} z_{N,k}(x) f(x)\, dx \ra\int_{\mathbb{R}^m}
z_{k}(x) f(x) \,dx = I_q , \qquad N\ra\infty.
\]
But this follows from Lebesgue's bounded convergence theorem, since
$z_{N,k}(x)$ is bounded (take $\delta=0$ in $J_N$).

\subsection{\texorpdfstring{Proof of Theorem \protect\ref{Th:consistency}}{Proof of Theorem 3.2}}

We shall use the same notations as in the proof of Theorem
\ref{Th:CV of mean} and write $\hat I_{N,k,q}=(1/N) \sum_{i=1}^N
\zeta_{N,i,k}^{1-q}$, so that
%
\begin{eqnarray}\label{main-consistency}
\Ex(\hat I_{N,k,q}-I_q)^2 &=& \frac{ \Ex
(\zeta_{N,i,k}^{1-q}-I_q)^2 }{N} \nonumber\\[-8pt]\\[-8pt]
&&{} + \frac{1}{N^2} \sum_{i\neq j} \Ex\{(\zeta_{N,i,k}^{1-q}-I_q)
(\zeta_{N,j,k}^{1-q}-I_q) \} . \nonumber
\end{eqnarray}
We consider the cases $q<1$ and $q>1$ separately.


$q<1$.\quad
Note that $2q-1<q<1$ and Lemma \ref{L:Iq}(ii) gives $I_q<\infty$
when $I_{2q-1}<\infty$. Consider the first term on the right-hand
side of (\ref{main-consistency}). We have
%
\begin{equation}\label{1st-term}
\Ex(\zeta_{N,i,k}^{1-q}-I_q)^2 = \Ex
(\zeta_{N,i,k}^{1-q})^2 +I_q^2-2I_q \Ex\zeta_{N,i,k}^{1-q},
\end{equation}
where the last term tends to $-2I_q^2$ from Theorem \ref{Th:CV of
mean}. Consider the first term,
\[
\Ex(\zeta_{N,i,k}^{1-q})^2 = \int_{\mathbb{R}^m} \int_0^\infty
u^{2(1-q)} f_{N,x,k}(u) f(x) \,du \,dx .
\]
Since the function $u\ra u^{1-q}$ is bounded on every bounded
interval, it tends to
\[
\int_{\mathbb{R}^m} \int_0^\infty u^{2(1-q)} f_{x,k}(u) f(x)
\,du \,dx = I_{2q-1} \frac{\Gamma(k+2-2q)\Gamma(k)}{\Gamma^2(k+1-q)}
\]
for any $q<(k+2)/2$ (generalized Helly--Bray lemma,
L\'oeve \cite{Loeve77}, page 187). Therefore, $\Ex
(\zeta_{N,i,k}^{1-q}-I_q)^2$ tends to a finite limit and the first
term on the right-hand side of (\ref{main-consistency}) tends to
zero as $N\ra\infty$.

Consider now the second term of (\ref{main-consistency}). We show
that
\begin{eqnarray*}
&&\Ex\{(\zeta_{N,i,k}^{1-q}-I_q) (\zeta_{N,j,k}^{1-q}-I_q) \}\\
&&\qquad =
\Ex\{\zeta_{N,i,k}^{1-q}\zeta_{N,j,k}^{1-q}\}+I_q^2-2I_q \Ex
\zeta_{N,i,k}^{1-q} \ra0 , \qquad N\ra\infty.
\end{eqnarray*}
Since $\Ex\zeta_{N,i,k}^{1-q}\ra I_q$ from Theorem \ref{Th:CV of
mean}, we only need to show that\break
$\Ex\{\zeta_{N,i,k}^{1-q}\zeta_{N,j,k}^{1-q}\} \ra I_q^2$. Define
\begin{eqnarray*}
F_{N,x,y,k}(u,v) &=& \Pr\{ \zeta_{N,i,k} < u , \ \zeta_{N,j,k}
< v | X_i=x,X_j=y\} , \\
&=& \Pr\bigl\{ \rho_{k,N-1}^{(i)} < R_N(u) , \ \rho_{k,N-1}^{(j)} < R_N(v)
| X_i=x,X_j=y\bigr\},
\end{eqnarray*}
so that
%
\begin{eqnarray}
\label{EE1}
&&\Ex\{\zeta_{N,i,k}^{1-q}\zeta_{N,j,k}^{1-q}\}\nonumber\\[-8pt]\\[-8pt]
&&\qquad= \int_{\mathbb{R}^m}
\int_{\mathbb{R}^m} \int_0^\infty\int_0^\infty u^{1-q} v^{1-q}
F_{N,x,y,k}(du,dv) f(x) f(y) \,dx \,dy .\nonumber
\end{eqnarray}
Let us assume that $x \neq y$. From the definition of
$R_N(u)$ [see (\ref{RN})] there exist $N_0=N_0(x,y,u,v)$ such that
$\SB[x,R_N(u)]\cap\SB[y,R_N(v)]=\varnothing$ for $N>N_0$ and
thus,
\begin{eqnarray*}
F_{N,x,y,k}(u,v) &=& \sum_{j=k}^{N-2} \sum_{l=k}^{N-2-j}
\pmatrix{ N-2 \cr j}
\pmatrix{ N-2-j \cr l}\\
&&\hspace*{47pt}{}\times
p_{N,u}^j p_{N,v}^l (1-p_{N,u}-p_{N,v})^{N-2-j-l}
\end{eqnarray*}
with $p_{N,u}=\int_{\SB[x,R_N(u)]} f(t)\,dt$,
$p_{N,v}=\int_{\SB[y,R_N(v)]} f(t)\,dt$. Hence, for $N>N_0$,
\begin{eqnarray*}
F_{N,x,y,k}(u,v) &=& F_{N-1,x,k}(u)+F_{N-1,y,k}(v)-1 \\
&& {} + \sum_{j=0}^{k-1} \sum_{l=0}^{k-1}
\pmatrix{ N-2 \cr j }
\pmatrix{ N-2-j \cr l }\\
&&\hspace*{44pt}{}\times
p_{N,u}^j p_{N,v}^l (1-p_{N,u}-p_{N,v})^{N-2-j-l} .
\end{eqnarray*}
Similarly to the proof of Theorem \ref{Th:CV of mean}, we then
obtain
%
\begin{equation}\label{EE2}
F_{N,x,y,k}(u,v) \ra F_{x,y,k}(u,v)=F_{x,k}(u)F_{y,k}(v) , \qquad
N\ra\infty,
\end{equation}
for $\mu_\SL$-almost any $x$ and $y$ with
%
\begin{equation}\label{EE3}
\int_0^\infty\int_0^\infty u^{1-q} v^{1-q}
F_{x,y,k}(du,dv)=f^{q-1}(x)f^{q-1}(y),
\end{equation}
for any $q<k+1$. Since the function $u\ra u^{1-q}$ is bounded on
every bounded interval, (\ref{EE1}) gives
\[
\Ex\{\zeta_{N,i,k}^{1-q}\zeta_{N,j,k}^{1-q}\} \ra
\int_{\mathbb{R}^m} \int_{\mathbb{R}^m} f^q(x) f^q(y) \,dx \,dy
= I_q^2 , \qquad N\ra\infty
\]
(generalized Helly--Bray lemma, \cite{Loeve77}, page 187). This
completes the proof that $\Ex(\hat I_{N,k,q}-I_q)^2 \ra0$.
Therefore, $\hat I_{N,k,q} \rap I_q$, when $N\ra\infty$.

$q>1$.\quad
Note that from Lemma \ref{L:Iq}(i) $I_q$ and $I_{2q-1}$ both exist.
Consider the first term on the right-hand side of
(\ref{main-consistency}). We have again (\ref{1st-term}) where the
last term tends to $-2I_q^2$ (the assumptions of the theorem imply
$q<k+1$ so that Theorem \ref{Th:CV of mean} applies). Consider the
first term of (\ref{1st-term}). Define
\[
J_N' = \int_0^\infty u^{2(1-q)(1+\delta)} F_{N,x,k}(du) ,
\]
we show that $\sup_N J_N'<\infty$ for some $\delta>0$. From the
assumptions of the theorem, $2q<k+2$. Let $\beta=2(1-q)(1+\delta)$,
so that $\beta<0$ and take $\delta<(k+2-2q)/[2(q-1)]$ so that
$\beta+k>0$. Using Lemma \ref{L:moment-F} and developments similar
to the proof of Theorem \ref{Th:CV of mean}, we obtain
\begin{eqnarray*}
J_N' &=& -\beta\int_0^\infty u^{\beta-1} F_{N,x,k}(du) \leq1-
\beta\int_0^1 u^{\beta-1} F_{N,x,k}(du) \\
& \leq & 1-\beta U_k \int_0^1 u^{k+\beta-1}\, du = 1 -
\frac{\beta U_k}{k+\beta} < \infty,
\end{eqnarray*}
where $U_k$ is given by (\ref{Uk}). Theorem 2.5.1 of
Bierens \cite{Bierens94} then implies
\begin{eqnarray*}
\int_0^\infty u^{2(1-q)} F_{N,x,k}(du) &\ra&\int_0^\infty
u^{2(1-q)} F_{x,k}(du) \\
&=& \frac{
\Gamma(k+2-2q)\Gamma(k)}{\Gamma^2(k+1-q)} f^{2q-2}(x)
\end{eqnarray*}
for $\mu$-almost any $x$, $q<(k+2)/2$ and Lebesgue's bounded
convergence theorem gives $\Ex(\zeta_{N,i,k}^{1-q})^2 \ra
I_{2q-1}\Gamma(k+2-2q)\Gamma(k)/\Gamma^2(k+1-q)$, $N\ra\infty$. The
first term of (\ref{main-consistency}) thus tends to zero.

Consider now the second term. As in the case $q<1$, we only need to
show that $\Ex\{\zeta_{N,i,k}^{1-q}\zeta_{N,j,k}^{1-q}\} \ra I_q^2$.
Define
\[
J''_N = \int_0^\infty\int_0^\infty u^{(1-q)(1+\delta)}
v^{(1-q)(1+\delta)} F_{N,x,y,k}(du,dv) .
\]
Using (\ref{EE2}, \ref{EE3}), proving that $\sup_N J''_N < J(x,y) <
\infty$ for some $\delta>0$ will then establish that
%
\begin{eqnarray}\label{EE4}
&&\int_0^\infty\int_0^\infty u^{1-q} v^{1-q} F_{N,x,y,k}(du,dv)\nonumber\\[-8pt]\\[-8pt]
&&\qquad\ra
f^{q-1}(x)f^{q-1}(y) , \qquad
N\ra\infty,\nonumber
\end{eqnarray}
for $\mu$-almost $x$ and $y$; see Theorem 2.5.1 of
Bierens \cite{Bierens94}. Using (\ref{EE1}), if
%
\begin{equation}\label{C1}
\int_{\mathbb{R}^m}\int_{\mathbb{R}^m} J(x,y) f(x) f(y) \,dx
\,dy<\infty,
\end{equation}
Lebesgue's dominated convergence theorem will then complete the
proof.

Integration by parts, as in the proof of Lemma \ref{L:moment-F},
gives
\[
J_N'' = \beta^2 \int_0^\infty\int_0^\infty u^{\beta-1}v^{\beta-1}
F_{N,x,y,k}(u,v) \,du \,dv ,
\]
where $\beta=(1-q)(1+\delta)<0$. We use different bounds for
$F_{N,x,y,k}(u,v)$ on three different parts of the $(u,v)$
plane.

(i) Suppose that $\max[R_N(u),R_N(v)]\leq\|x-y\|$, which is
equivalent to $(u,v)\in\SD_1=[0,\Lambda]\times[0,\Lambda]$ with
$\Lambda=\Lambda(k,N,x,y)=(N-1)V_mC_k \|x-y\|^m$. This means that
the balls $\SB[x,R_N(u)]$ and $\SB[y,R_N(v)]$ either do not
intersect, or, when they do, their intersection contains neither $x$
nor $y$. In that case, we use
\[
F_{N,x,y,k}(u,v) < \min[F_{N-1,x,k}(u),F_{N-1,y,k}(v)] <
F^{1/2}_{N-1,x,k}(u) F^{1/2}_{N-1,y,k}(v)
\]
and
\begin{eqnarray*}
J_N''^{(1)} &=& \beta^2 \int_{\SD_1} u^{\beta-1}v^{\beta-1}
F_{N,x,y,k}(u,v) \,du \,dv \\
& < &\beta^2 \biggl[\int_0^\Lambda u^{\beta-1}
F^{1/2}_{N-1,x,k}(u)\, du\biggr] \biggl[\int_0^\Lambda v^{\beta-1}
F^{1/2}_{N-1,y,k}(v)\, dv\biggr] \\
& <& \beta^2 \biggl[ U_k^{1/2} \int_0^1 u^{\beta-1+k/2} \,du +
\int_1^\infty u^{\beta-1} \,du \biggr]^2 \\
& =& \beta^2 \biggl[ U_k^{1/2} \frac{2}{2\beta+k} - \frac{1}{\beta}
\biggr]^2 < \infty,
\end{eqnarray*}
where we used the bound (\ref{Uk}) for $F_{N-1,x,k}(u)$ when $u<1$,
$F_{N-1,x,k}(u)<1$ for $u\geq1$ and choose
$\delta<(k+2-2q)/[2(q-1)]$ so that $2\beta+k>0$ [this choice of
$\delta$ is legitimate since $q<(k+2)/2$].

(ii) Suppose, without any loss of generality, that $u<v$ and
consider the domain defined by $R_N(u) \leq\|x-y\|<R_N(v)$, that
is, $(u,v)\in\SD_2=[0,\Lambda]\times(\Lambda,\infty)$. The cases
$k=1$ and $k\geq2$ must be treated separately since $\SB[y,R_N(v)]$
contains~$x$.

When $k=1$, $F_{N,x,y,1}(u,v) = F_{N-1,x,1}(u)$ and we have
%
\begin{eqnarray} \label{J2}
J_N''^{(2)} &=& \beta^2 \int_{\SD_2} u^{\beta-1}v^{\beta-1}
F_{N,x,y,1}(u,v)\, du \,dv \nonumber\\
& < & \beta^2 \biggl[\int_0^\Lambda u^{\beta-1}
F_{N-1,x,1}(u) \,du\biggr] \biggl[\int_\Lambda^\infty v^{\beta-1}
\,dv\biggr] \nonumber\\
& <& \beta^2 \biggl[ U_1 \int_0^1 u^{\beta} du +
\int_1^\infty u^{\beta-1} \,du \biggr] \biggl(-\frac{\Lambda^\beta
}{\beta}\biggr) \\
& = &-\beta\biggl[ \frac{U_1}{\beta+1} - \frac{1}{\beta}
\biggr] \Lambda^\beta\nonumber\\
& <& J^{(2)}(x,y) = -\beta\biggl[ \frac{U_1}{\beta+1} - \frac{1}{\beta}
\biggr] V_m^\beta C_1^\beta\|x-y\|^{m\beta} ,\nonumber
\end{eqnarray}
where we used (\ref{Uk}) and take $\delta<(2-q)/(q-1)$ so that
$\beta>-1$ (this choice of $\delta$ is legitimate since $q<2$).

Suppose now that $k\geq2$. We have $F_{N,x,y,k}(u,v) <
F^{1-\alpha}_{N-1,x,k}(u) F^{\alpha}_{N-1,y,k-1}(v)$, $\forall
\ma\in(0,1)$. Developments similar to those used for the derivation
of (\ref{Uk}) give for $v<1$
%
\begin{eqnarray}\label{Vk-1}
&&\frac{F_{N-1,y,k-1}(v)}{v^{k-1}}\nonumber\\[-8pt]\\[-8pt]
 &&\qquad< V_{k-1} = \frac{\bar
f^{k-1}}{C_k^{k-1} (k-1)!}
 + \frac{\bar f^{k-1}}{C_k^{k-1}}
\biggl\{\exp\biggl[ \frac{\bar f}{C_k} \biggr]-1\biggr\} .\nonumber
\end{eqnarray}
We obtain
\begin{eqnarray*}
J_N''^{(2)} &=& \beta^2 \int_{\SD_2} u^{\beta-1}v^{\beta-1}
F_{N,x,y,k}(u,v) \,du \,dv \\
&< & \beta^2 \biggl[\int_0^\Lambda u^{\beta-1}
F^{1-\ma}_{N-1,x,k}(u) \,du\biggr] \biggl[\int_\Lambda^\infty v^{\beta-1}
F^{\ma}_{N-1,y,k-1}(v) \,dv\biggr] \\
& < &\beta^2 \biggl[ U_k^{1-\ma} \int_0^1 u^{\beta-1+(1-\ma)k} \,du +
\int_1^\infty u^{\beta-1} \,du \biggr] \\
&& {} \times\biggl[ V_{k-1}^\ma\int_0^1 v^{\beta-1+(k-1)\ma
} \,dv +
\int_1^\infty v^{\beta-1}\, dv \biggr] \\
& = & \beta^2 \biggl[ \frac{U_k^{1-\ma}}{k(1-\ma)+\beta} - \frac
{1}{\beta}
\biggr] \biggl[ \frac{V_{k-1}^\ma}{(k-1)\ma+\beta} - \frac{1}{\beta}
\biggr] < \infty,
\end{eqnarray*}
where we used (\ref{Uk}, \ref{Vk-1}) and require $\beta+k(1-\ma)>0$
and $\beta+(k-1)\ma>0$. For that we take $\ma=\ma_k=k/(2k-1)$.
Indeed, from the assumptions of the theorem,
$q<(k+1)/2<(k^2+k-1)/(2k-1)$ so that we can choose
$\delta<[(k^2+k-1)-q(2k-1)]/[(q-1)(2k-1)]$, which ensures that both
$\beta+k(1-\ma_k)>0$ and $\beta+(k-1)\ma_k>0$.

(iii) Suppose finally that $\|x-y\|< \min[R_N(u),R_N(v)]$, that is,
$(u,v)\in\SD_3=(\Lambda,\infty)\times(\Lambda,\infty)$. In that
case, each of the balls $\SB[x,R_N(u)]$ and $\SB[y,R_N(v)]$ contains
both $x$ and $y$. Again, the case $k=1$ and $k\geq2$ must be
distinguished.

When $k=1$, $F_{N,x,y,1}(u,v)=1$ and
%
\begin{eqnarray}\label{J3}
J_N''^{(3)} &=& \beta^2 \int_{\SD_3} u^{\beta-1}v^{\beta-1}
F_{N,x,y,1}(u,v) \,du\, dv \nonumber\\
&=& \beta^2 \biggl[ \int_\Lambda^\infty u^{\beta-1} \,du \biggr]^2
=
\Lambda^{2\beta} \\
&<& J^{(3)}(x,y)= V_m^{2\beta} C_1^{2\beta}
\|x-y\|^{2m\beta} .\nonumber
\end{eqnarray}

When $k\geq2$, $F_{N,x,y,k}(u,v) < F^{1/2}_{N-1,x,k-1}(u)
F^{1/2}_{N-1,y,k-1}(v)$ and
\begin{eqnarray*}
J_N''^{(3)} &=& \beta^2 \int_{\SD_3} u^{\beta-1}v^{\beta-1}
F_{N,x,y,k}(u,v) \,du \,dv \\[0.5pt]
&< & \beta^2 \biggl[\int_\Lambda^\infty u^{\beta-1}
F^{1/2}_{N-1,x,k-1}(u) \,du\biggr] \\[0.5pt]
&&{}\times\biggl[\int_\Lambda^\infty v^{\beta-1}
F^{1/2}_{N-1,y,k-1}(v) \,dv\biggr] \\[0.5pt]
&< & \beta^2 \biggl[ V_{k-1}^{1/2} \frac{2}{2\beta+k-1} - \frac{1}{\beta}
\biggr]^2 < \infty,
\end{eqnarray*}
where we used (\ref{Vk-1}) and take $\delta<[(k+1)-2q]/[2(q-1)]$ so
that $k-1+2\beta>0$ [this choice of $\delta$ is legitimate since
$q<(k+1)/2$].

Summarizing the three cases above, we obtain
$J_N''=J_N''^{(1)}+2J_N''^{(2)}+J_N''^{(3)}$ with different bounds
for $J_N''^{(2)}$ and $J_N''^{(3)}$ depending on whether $k=1$ or
$k\geq2$. This proves (\ref{EE4}).

When $k\geq2$, the bound on $J_N''$ does not depend on $x,y$ and
Lebesgue's bounded convergence theorem implies
$\Ex\{\zeta_{N,i,k}^{1-q}\zeta_{N,j,k}^{1-q}\} \ra I_q^2$, which
completes the proof of the theorem; see (\ref{C1}).

When $k=1$, the condition (\ref{C1}) is satisfied if $2\beta>-1$
[see (\ref{J2}), (\ref{J3})], which is ensured by the choice
$\delta<(3-2q)/[2(q-1)]$ (legitimate since $q<3/2$). Indeed, we
can write
\[
\int_{\mathbb{R}^m} \int_{\mathbb{R}^m} \|x-y\|^{\mg} f(x)
f(y) \,dx \,dy = \int_{\mathbb{R}^m} \|x\|^{\mg} g(x) \,dx ,
\]
where $g(x) = \int_{\mathbb{R}^m} f(x+y) f(y)\, dy$, and thus
(since $\mg<0$),
\begin{eqnarray*}
\int_{\mathbb{R}^m} \int_{\mathbb{R}^m} \|x-y\|^\mg f(x)
f(y) \,dx \,dy &\leq& {\bar f}^2 \int_{\|x\|<1} \|x\|^{\mg}
dx + I_2 \\
& = & {\bar f}^2 \frac{mV_m}{\mg+m} + I_2,
\end{eqnarray*}
when $\mg>-m$. When $\delta<(3-2q)/[2(q-1)]$, Lebesgue's dominated
convergence theorem thus implies
$\Ex\{\zeta_{N,i,k}^{1-q}\zeta_{N,j,k}^{1-q}\} \ra I_q^2$ , which
completes the proof of the theorem.


\subsection{\texorpdfstring{Proof of Corollary \protect\ref{C:biasH1}}{Proof of Corollary 3.2}}

The existence of $H_1$ directly follows from that of $I_{q_1}$ for
$q_1<1$ and the boundedness of $f$.

\subsubsection*{Asymptotic unbiasedness} We have
\[
\Ex\hat H_{N,k,1}=\Ex\log\xi_{N,i,k}
=\Ex[\Ex(\log\xi_{N,i,k}|X_i=x)] ,
\]
where the only difference between the random variables
$\zeta_{N,i,k}$ (\ref{xi}) and $\xi_{N,i,k}$ (\ref{zeta}) is the
substitution of $\exp[-\Psi(k)]$ for $C_k$. Similarly to the proof
of Theorem~\ref{Th:CV of mean}, we define $ F_{N,x,k}(u) =
\Pr(\xi_{N,i,k} < u|X_i=x)
= \Pr[\rho_{k,N-1}^{(i)} < R_N(u)|X_i=x]$
with now $R_{N}(u)=(u/\{(N-1)V_m \exp[-\Psi(k)]\})^{1/m}$. Following
the same steps as in the proof of Theorem \ref{Th:CV of mean}, we
then obtain
\[
F_{N,x,k}(u) \ra F_{x,k}(u)=1-\exp(-\ml u) \sum_{j=0}^{k-1}
\frac{(\ml u)^j}{j!} , \qquad N\ra\infty,
\]
for $\mu_\SL$-almost any $x$, with $\ml=f(x)\exp[\Psi(k)]$.

Direct calculation gives $\int_0^\infty\log(u) F_{x,k}(du) = - \log
f(x)$. We shall use again Theorem 2.5.1 of Bierens \cite{Bierens94}, page
34, and show that
%
\begin{equation}\label{JN}
J_N = \int_0^\infty|\log(u)|^{1+\delta} F_{N,x,k}(du) < \infty
,
\end{equation}
for some $\delta>0$, which implies
\[
\int_0^\infty\log(u) F_{N,x,k}(du) \ra\int_0^\infty\log(u)
F_{x,k}(du) = - \log f(x) , \qquad N\ra\infty,
\]
for $\mu_\SL$-almost any $x$. The convergence
\[
\int_{\mathbb{R}^m} \int_0^\infty\log(u) F_{N,x,k}(du) f(x)
\,dx \ra H_1 , \qquad N\ra\infty,
\]
then follows from Lebesgue's bounded convergence theorem.

In order to prove (\ref{JN}), we write
%
\begin{equation}\label{JN2}
J_N = \int_0^1 |\log(u)|^{1+\delta} F_{N,x,k}(du) +
\int_1^\infty|\log(u)|^{1+\delta} F_{N,x,k}(du) .
\end{equation}
Since $f$ is bounded, we can take $q_2>1$ (and smaller than $k+1$)
such that $\int_0^\infty u^{1-q_2} F_{N,x,k}(du) < \infty$; see
(\ref{resT1}). Since $|\log(u)|^{1+\delta}/u^{1-q_2}\ra0$ when
$u\ra0$, it implies that the first integral on the right-hand side
of (\ref{JN2}) is finite. Similarly, since, by assumption, $I_{q_1}$
exists for some $q_1<1$, $\int_0^\infty u^{1-q_1} F_{N,x,k}(du) <
\infty$ and $|\log(u)|^{1+\delta}/u^{1-q_1}\ra0$, $u\ra\infty$,
implies that the second integral on the right-hand side of
(\ref{JN2}) is finite, which completes the proof that $\Ex\hat
H_{N,k,1} \ra H_1$ as $N\ra\infty$.\looseness=1

\subsubsection*{$L_2$ consistency}
Similarly to the proof of asymptotic
unbiasedness, we only need
to replace $\zeta_{N,i,k}$ (\ref{xi}) by $\xi_{N,i,k}$ (\ref{zeta})
and $C_k$ by $\exp[-\Psi(k)]$ in the proof of Theorem
\ref{Th:consistency}. When we now compute
%
\begin{eqnarray} \label{main-consistency-Shannon}
\quad \Ex( \hat H_{N,k,1} - H_1)^2 &=& \frac{ \Ex( \log\xi_{N,i,k} -
H_1)^2}{N} \nonumber\\[-8pt]\\[-8pt]
&&{} + \frac{1}{N^2} \sum_{i\neq j} \Ex\{ ( \log\xi_{N,i,k} - H_1) (
\log\xi_{N,j,k} - H_1) \} , \nonumber
\end{eqnarray}
in the first term, $\Ex( \log\xi_{N,i,k} - H_1)^2$ tends to
\[
\int_{\mathbb{R}^m} \log^2 f(x) f(x)\, dx - H_1^2 + \dot{\Psi}(k)
= \var[\log f(X)] + \dot{\Psi}(k) ,
\]
where $\dot{\Psi}(z)$ is the trigamma function,
$\dot{\Psi}(z)=d^2\log\Gamma(z)/dz^2$, and for the second term the
developments are similar to those in Theorem \ref{Th:consistency}.
For instance, equation (\ref{EE4}) now becomes $ \int_0^\infty
\int_0^\infty\log u \log v F_{N,x,y,k}(du,dv) \ra\log f(x)
\log f(y) $, $N\ra\infty$,
for $\mu$-almost $x$ and $y$. We can then show that $\Ex\{\log
\xi_{N,i,k} \log\xi_{N,j,k}\} \ra H_1^2$, so that $\Ex( \hat
H_{N,k,1} - H_1)^2 \ra0$, $N\ra\infty$.

\section*{Acknowledgments}
The authors wish to
thank Anatoly A.\ Zhigljavsky from the Cardiff School of Mathematics
for helpful discussions. Comments from the associate editor and
the referees helped in improving the presentation of our
results.

\printaddresses

\end{document}